\documentclass[MSNbibl,number,citesort,seceqn,dvips]{arxbj}
\usepackage{graphicx}
\aid{0}
\volume{19}
\issue{5A}
\pubyear{2013}
\firstpage{1688}
\lastpage{1713}
\doi{10.3150/12-BEJ426} 

\makeatletter
\newtheorem{lemma}{Lemma}[section]
\newremark{example}[lemma]{Example}
\newtheorem{theorem}[lemma]{Theorem}
\newremark{remark}[lemma]{Remark}

\newcommand{\eqref}[1]{(\ref{#1})}
\newcommand{\eps}{\varepsilon}
\newcommand{\ov}{\overline}
\newcommand{\wt}{\widetilde}
\newcommand{\garch}{$\operatorname{GARCH}(1,1)$}
\newcommand{\bfR}{{\mathbf{R}}}
\newcommand{\sign}{\operatorname{sign}}
\newcommand{\bfTh}{{\boldsymbol{\Theta}}}
\newcommand{\std}{\stackrel{d}{\rightarrow}}
\newcommand{\stv}{\stackrel{v}{\rightarrow}}
\newcommand{\stw}{\stackrel{w}{\rightarrow}}
\newcommand{\eqd}{\stackrel{d}{=}}
\newcommand{\nto}{n\to\infty}
\newcommand{\xto}{x\to\infty}
\newcommand{\wh}{\widehat}
\newcommand{\vep}{\varepsilon}
\newcommand{\la}{\lambda}
\newcommand{\bbr}{{\mathbb R}}
\newcommand{\bbz}{{\mathbb Z}}
\newcommand{\bbs}{{\mathbb S}}
\newcommand{\bfX}{{\mathbf X}}
\newcommand{\bfB}{{\mathbf B}}
\newcommand{\bfY}{{\mathbf Y}}
\newcommand{\bfA}{{\mathbf A}}
\newcommand{\bfe}{{\mathbf e}}
\newcommand{\bfc}{{\mathbf c}}
\newcommand{\bfs}{{\mathbf s}}
\makeatother

\begin{document}
\begin{frontmatter}

\title{Stochastic volatility models with possible extremal clustering}
\runtitle{Stochastic volatility models with extremal clustering}

\begin{aug}
\author[1]{\fnms{Thomas} \snm{Mikosch}\corref{}\thanksref{1}\ead[label=e1]{mikosch@math.ku.dk}} \and
\author[2]{\fnms{Mohsen} \snm{Rezapour}\thanksref{2}\ead[label=e2]{mohsenrzp@gmail.com}}
\runauthor{T. Mikosch and M. Rezapour} 
\address[1]{University of Copenhagen, Department of Mathematics,
Universitetsparken 5,
DK-2100 Copenhagen, Denmark. \printead{e1}}
\address[2]{Department of Statistics, Shahid Bahonar University,
Kerman, Iran.\\ \printead{e2}}
\end{aug}

\received{\smonth{5} \syear{2011}}
\revised{\smonth{12} \syear{2011}}

%
\begin{abstract}
In this paper we consider a heavy-tailed stochastic volatility model,
$X_t=\sigma_tZ_t$,
$t\in\mathbb{Z}$,
where the volatility sequence $(\sigma_t)$
and the i.i.d. noise sequence $(Z_t)$ are
assumed independent, $(\sigma_t)$ is regularly varying with index
$\alpha>0$,
and the $Z_t$'s have moments of order larger than $\alpha$. In
the literature (see \textit{Ann. Appl. Probab.} \textbf{8} (1998)
664--675, \textit{J. Appl.
Probab.} \textbf{38A} (2001) 93--104, In \textit{Handbook of
Financial Time Series} (2009)
355--364 Springer), it is typically
assumed that $(\log\sigma_t)$ is a Gaussian stationary sequence and
the $Z_t$'s are regularly varying with some index $\alpha$ (i.e.,
$(\sigma_t)$
has lighter tails than the $Z_t$'s), or that $(Z_t)$ is i.i.d. centered
Gaussian. In these cases, we see that the sequence $(X_t)$ does not
exhibit extremal clustering. In contrast to this situation, under the
conditions of this paper, both situations are possible; $(X_t)$
may or may not have extremal clustering, depending on the clustering
behavior of the $\sigma$-sequence.
\end{abstract}

%

\end{frontmatter}

\section{Introduction}
The \emph{stochastic volatility model}
%
\begin{eqnarray}\label{eq:sv}
X_t=\sigma_t Z_t , \qquad t\in\bbz,
\end{eqnarray}
has attracted some attention in the financial
time series literature. Here the \emph{volatility sequence} $(\sigma
_t)$ is (strictly)
stationary and consists of
non-negative random variables independent of the i.i.d. sequence $(Z_t)$. We refer to
\cite{andersen:davis:kreiss:mikosch:2009} for a recent overview of the
theory of stochastic volatility models. The popular GARCH model has
the same
structure \eqref{eq:sv}, but every $Z_t$ feeds into the future
volatilities $\sigma_{t+k}$, $k\ge1$, and thus
$(\sigma_t)$ and $(Z_t)$ are dependent in this case; see, for example, the
definition of a \garch\ process in Example~\ref{exam:garch}. However, neither
$\sigma_t$ nor $Z_t$ is directly observable, and thus whether we
prefer a stochastic volatility a GARCH, or any
other model for returns depends on our modeling efforts.

Previous research on extremes (e.g., \cite
{breidt:davis:1998,davis:mikosch:2001,davis:mikosch:2009,kulik:soulier:2011})
has focused mainly on stochastic volatility models, where $(\log
\sigma_t)$ constitutes a Gaussian stationary process and $(Z_t)$ is
light-tailed (e.g., centered Gaussian) or rather heavy-tailed in the
sense that there exists $\alpha>0$, a slowly varying function $L$
and constants
$p,q\ge0$ such that $p+q=1$ and
%
\begin{eqnarray}\label{eq:rv}
P(Z>x)\sim p x^{-\alpha}L(x) \quad\mbox{and}\quad
P(Z\le-x)\sim q x^{-\alpha}L(x),\qquad \xto.
\end{eqnarray}
Here and in what follows, for any (strictly) stationary sequence $(Y_t)$,
$Y$ denotes a generic element. A random variable $Z$ satisfying \eqref{eq:rv}
will be called \emph{regularly varying with index}~$\alpha$.

Under the foregoing conditions, the sequence $(X_t)$ exhibits extremal behavior
similar to an i.i.d. sequence whatever the strength of dependence in
the Gaussian
log-volatility sequence. In particular, $(X_t)$ does not have extremal
clusters. It is common to measure extremal clustering in a stationary
sequence $(Y_t)$ by considering the \emph{extremal index}; suppose that
an i.i.d. sequence $(\wt Y_t)$ with the same marginal distribution as $Y$
satisfies the limit relation
\[
\lim_{\nto} P\bigl(c_n^{-1}\bigl(\max(\wt Y_1,\ldots,\wt Y_n)-d_n\bigr)\le
x\bigr)=H(x),\qquad x\in\bbr
\]
for suitable constants $c_n>0$, $d_n\in\bbr$ and a nondegenerate
limit distribution function $H$ (which is necessarily continuous). If
the same limit
relation holds with $\max(\wt Y_1,\ldots,\wt Y_n)$ replaced by
$\max(Y_1,\ldots, Y_n)$ and $H$ replaced by
$H^\theta$ for some $\theta\in[0,1]$, then
$\theta$ is called the \emph{extremal index} of $(Y_t)$. Clearly,
that the smaller the $\theta$, the stronger the extremal clustering effect
present in the sequence. Under the aforementioned conditions, the
stochastic volatility model
$(X_t)$ has extremal index 1; that is, this process does not exhibit
extremal clustering. However, real-life financial returns typically cluster
around large positive and small negative values. This effect is
described by
the GARCH model, which under general conditions has an extremal index
$\theta<1$ (see \cite{basrak:davis:mikosch:2002,mikosch:starica:2000}).

The aim of this paper is to show that the lack of extremal clustering
in stochastic volatility models is due to the conditions on the tails
of distributions of the
sequences $(\sigma_t)$ and $(Z_t)$. In particular, we focus on the
heavy-tailed situation when the distribution of $\sigma$ has power
law tails
in the sense that there exist $\alpha>0$ and a slowly varying function $L$
such that
\[
P(\sigma>x) \sim x^{-\alpha} L(x),\qquad \xto;
\]
that is, $\sigma$ is regularly varying with index $\alpha$,
and $Z$ has lighter tail in the sense that
$E|Z|^{\alpha+\vep}<\infty$ for some $\vep>0$. By a result of Breiman
\cite{breiman:1965}, we then have
\[
P(X>x)\sim EZ_+^\alpha P(\sigma>x) \quad\mbox{and}\quad
P(X\le-x)\sim EZ_-^\alpha P(\sigma>x),\qquad \xto.
\]
This means that the tail behavior of $X$ is essentially determined by
the right tail of $\sigma$. This is in contrast to the situation
mentioned earlier. In that case, also by Breiman's result,
$P(X>x)\sim E\sigma^\alpha P(Z>x)$. The latter relation is
responsible for the lack of clustering; it indicates that extreme
values of the sequence $(X_t)$ are essentially\vadjust{\goodbreak} determined by the extremes
in the i.i.d. sequence $(Z_t)$, an extremal index $\theta=1$ result.
We mention in passing that extremal clustering also can be expected when
both the tails of $Z$ and $\sigma$ are regularly varying with the
same index
$\alpha>0$. In that case it is well known (see
\cite{embrechts:veraverbeke:1982}) that $X$ has regularly varying tails with a
slowly varying function $L$, which in general is rather difficult to
determine. We will not treat this case because it is of limited
interest
and will lead to
rather technical conditions.

The paper is organized as follows. In Section~\ref{sec:1} we
introduce the notion of a regularly varying sequence and review point
process convergence for
such a sequence which was developed by~\cite{davis:hsing:1995}. We then
state a result (Theorem~\ref{thm:main}) which translates mixing and
regular variation of the sequence $(\sigma_t)$ to the stochastic
volatility model $(X_t)$ defined
in \eqref{eq:sv}. Our results in
Sections~\ref{subsec:ar}--\ref{subsec:3} are concerned with
three major examples.
In Section~\ref{subsec:ar} we study the stochastic volatility model
\eqref{eq:sv}, where $(\sigma_t)$ is an exponential $\operatorname{AR}(1)$ process with
regularly varying marginals. We show that this model does not exhibit extremal
clustering, due to the lack of extremal clustering in $(\sigma_t)$.
We also show that an EGARCH model with the same volatility
dynamics has no extremal clustering either.
In Section~\ref{subsec:2} we assume that a positive power of
$(\sigma_t)$ satisfies a random coefficient autoregressive equation,
which we call
stochastic recurrence equation. In this case, the extremal clustering
of $(\sigma_t)$
translates to the stochastic volatility model. In Section~\ref
{subsec:3} we consider
another model with genuine extremal clustering. Here we assume that
$(\sigma_t)$ is some positive power of the absolute values of a
regularly varying moving average process.

\section{Preliminaries}\label{sec:1}
\subsection{Regularly varying sequences}
A strictly stationary sequence $(X_t)$ is said to be regularly
varying with
index $\alpha>0$ if for every $d\ge1$,
the vector $\bfX_d=(X_1,\ldots,X_d)'$ is regularly varying with index
$\alpha>0$. This means that there exists a sequence $(a_n)$ with
$a_n\to\infty$ and a sequence of non-null Radon measures $(\mu_d)$
on the Borel
$\sigma$-field of $\ov\bbr{}^d_0=\ov\bbr{}^d\setminus\{0\}$
such that for every $d\ge1$,
\[
n P(a_n^{-1}\bfX_d\in\cdot)\stv\mu_d(\cdot) ,
\]
where $\stv$ denotes vague convergence and $\mu_d$ satisfies the scaling
property
$\mu_d(t\cdot)=t^{-\alpha}\mu_d(\cdot)$, $t>0$.
The latter property justifies the term ``regular variation with index
$\alpha>0$.''
The sequence $(a_n)$ can be chosen as
such that $nP(|X_1|>a_n)\to1$. We refer to \cite
{resnick:1987,resnick:2007} for more
reading on regular variation and vague convergence of measures.
Examples of regularly varying sequences are GARCH
processes with i.i.d. Student or normal noise and ARMA processes with
i.i.d. regularly varying noise. \cite{davis:mikosch:2001} studied
the extremes of the stochastic volatility model
\eqref{eq:sv} under the assumptions that
$E\sigma^{\alpha+\vep}<\infty$ for some $\vep>0$
and $(Z_t)$ is i.i.d. regularly varying with index
$\alpha>0$. Then $(X_t)$ is regularly varying with index $\alpha$,
and the
measures $\mu_d$ are concentrated on the axes. This property is shared
with an i.i.d. regularly varying sequence $(X_t)$.

In this paper, we consider the opposite situation. We assume that
$(\sigma_t)$ is regularly varying with index $\alpha>0$, normalizing
constants
$(a_n)$ such that $n P(\sigma>a_n)\to1$, and limiting measures
$\nu_d$, $d=1,2,\ldots$\,. This means that for
${\boldsymbol\Sigma}_d=(\sigma_1,\ldots,\sigma_d)'$, $d\ge1$, the relations
\[
n P(a_n^{-1}{\boldsymbol\Sigma}_d\in\cdot)\stv\nu_d(\cdot)
\]
hold.
We also assume that
$E|Z|^{\alpha+\vep}<\infty$ for some $\vep>0$.
\begin{lemma}\label{lem:1}
Under the foregoing conditions, $(X_t)$ is regularly varying with index
$\alpha$
and limiting measures $\mu_d$, $d=1,2,\ldots$\,, given by the relation
%
\begin{equation}\label{eq:ms}
\mu_d(\cdot)= E\nu_d\{\bfs\in\bbr^d_+\dvt (Z_1s_1,\ldots,Z_d
s_d)\in\cdot\} .
\end{equation}
\end{lemma}
\begin{pf}
Assuming that all vectors are written in column form, we have
$
\bfX_d=\bfA{\boldsymbol
\Sigma}_d
$,
where $\bfA=\operatorname{diag}(Z_1,\ldots,Z_d)$.
The matrix $\bfA$ has moment of order $\alpha+\vep$ and then regular
variation
 of $\bfX_d$ with normalizing constants
$(a_n)$ given by $n P(\sigma>a_n)\to1$
and the form of the limit measures $\mu_d$ follow from the
multivariate Breiman result (see \cite{basrak:davis:mikosch:2002}).
\end{pf}

The limits \eqref{eq:ms} are generally difficult to evaluate. We
consider some simple examples.

\begin{example}\label{exam:limit}
Assume that $\nu_d$ is concentrated
on the axes, that is, it has
the form
\[
\nu_d(\cdot)=c_d \sum_{i=1}^d \int_0^\infty x^{-\alpha-1} I_{\{
x\bfe_i\in\cdot\}}\,\mathrm{d}x
\]
for some constants $c_d>0$, where $ \bfe_i$ denotes the $i$th unit
vector in $\bbr^d$. Then \eqref{eq:ms} reads as
\begin{eqnarray*}
\mu_d(\cdot)&=& c_d\sum_{i=1}^d \int_0^\infty x^{-\alpha-1} P(x
Z_1\bfe_i\in\cdot) \,\mathrm{d}x
\\
&=&c_d\Biggl [EZ_+^\alpha\sum_{i=1}^d \int_0^\infty
x^{-\alpha-1} I_{\{x \bfe_i\in\cdot\}} \,\mathrm{d}x
+EZ_-^\alpha\sum_{i=1}^d \int_0^\infty
x^{-\alpha-1} I_{\{-x \bfe_i\in\cdot\}} \,\mathrm{d}x \Biggr] .
\end{eqnarray*}
\end{example}

Sometimes it is possible to characterize the limit measures $\mu_d$ by
their values on all sets of the form $A_\bfc=\{{\mathbf x}\in\bfX^d\dvt
\bfc'{\mathbf x}
>1\}$ for any choice of $\bfc$ in the unit sphere $\bbs^{d-1}$ of
$\bbr$ with respect to the Euclidean norm. However, in general, $\mu
_d$ cannot be reconstructed from its values on the sets $A_\bfc$ (see
\cite{basrak:davis:mikosch:2002a,boman:lindskog:2007,hult:lindskog:2006}).
\begin{example}
Consider an i.i.d. sequence of symmetric
$\beta$-stable
random variables
$(Z_t)$; that is, the characteristic function of $Z$ is given by
$\mathrm{e}^{-c |z|^{\beta}}$,
$z\in\bbr$,
for some $c>0$.
Assume that $\beta=2$ for $\alpha\ge2$ and
$2\ge\beta>\alpha$ for $\alpha<2$.
Then, for $\bfc\in\bbs^{d-1}$,
\begin{eqnarray*}
\mu_d(A_\bfc)&=& E\nu_d\Biggl\{\bfs\in\bbr^d_+\dvt \sum_{i=1}^d
c_iZ_is_i>1\Biggr\}
\\
&=&E\nu_d\Biggl\{\bfs\in\bbr^d_+\dvt Z
\Biggl(\sum_{i=1}^d |c_i|^\beta s_i^\beta\Biggr)^{1/\beta}>1\Biggr\}
\\
&=&EZ_+^\alpha\nu_d\Biggl\{\bfs\in\bbr^d_+\dvt
\Biggl(\sum_{i=1}^d |c_i|^\beta s_i^\beta\Biggr)^{1/\beta}>1\Biggr\}.
\end{eqnarray*}
The measure $\mu_d$, for example, is uniquely determined by the values
$\mu_d(A_\bfc)$, $\bfc\in\bbs^{d-1}$, provided that they are
positive and
$\alpha$ is not an integer (see \cite{basrak:davis:mikosch:2002a})
or, in view of the symmetry of the underlying distributions, if
$\alpha$ is an odd integer (see \cite{kluppelberg:pergam:2007}).
By virtue of the foregoing calculations, this means that $\mu_d$ is
uniquely determined
by the values of $\nu_d$ on the sets
%
\begin{eqnarray}\label{eq:el}
\Biggl\{\bfs\in\bbr^d_+\dvt
\Biggl(\sum_{i=1}^d |c_i|^\beta s_i^\beta\Biggr)^{1/\beta}>1\Biggr\}, \qquad
\bfc\in\bbs^{d-1},
\end{eqnarray}
provided that these values are positive.
For $\beta=2$, $Z$ is centered normal, and then \eqref{eq:el}
describes the complements of all ellipsoids in $\bbr^d$ with
$\sum_{i=1}^dc_i^2=1$ intersected with $\bbr_+^d$.
\end{example}
%
\subsection{Mixing conditions}\label{subsec:mix}
For the reader's convenience, here we introduce mixing concepts
for strictly stationary sequences $(X_t)$
used in this work. For $h\ge1$, let
\begin{eqnarray*}
\alpha_h&=&\sup_{A\in\sigma_{(-\infty,0]},B\in
\sigma_{[h,\infty)}}|P(A\cap B)-P(A)P(B)| ,
\\
\beta_h&=&E \Bigl(\sup_{B\in
\sigma_{[h,\infty)}}\bigl|P\bigl(B\vert\sigma_{(-\infty,0]}\bigr)-P(B)\bigr|\Bigr ),
\end{eqnarray*}
where $\sigma_A$ is the $\sigma$-field generated by $(X_t)_{t\in A}$ for
any $A\subset\bbz$.
The sequence $(X_t)$ is \emph{strongly mixing with rate function $(\alpha
_h)$} if
$\alpha_h\to0$ as $h\to\infty$. If $\beta_h\to0$ as $h\to\infty
$, then
$(X_t)$ is \emph{$\beta$-mixing with rate function $(\beta_h)$}.
Strong mixing is known to imply $\beta$-mixing (see Doukhan
\cite{doukhan:1994} for examples and comparisons of different mixing concepts).

Strong mixing and $\beta$-mixing were introduced in the context of
the central limit theory for partial sums of $(X_t)$. For partial
maxima of $(X_t)$,
other mixing concepts are more suitable (see, e.g., the conditions $D$
and $D'$
in
Leadbetter \textit{et al.} \cite{leadbetter:lindgren:rootzen:1983}). In this
paper, we make use of the condition~${\mathcal A}(a_n)$ introduced by
Davis and Hsing \cite{davis:hsing:1995}: Assume that there exists a
sequence $r_n\to\infty$ such that $r_n=\mathrm{o}(n)$ and
%
\begin{equation}\label{eq:aofa}
\Psi_f(N_n)- (\Psi_f(N_{n,r_n}))^{n/r_n}\to0 ,
\end{equation}
where $N_n$ is the point process of the points
$(a_n^{-1}X_t)_{t=1,\ldots,n}$, $N_{n,r_n}$ is the point process of
the points
$(a_n^{-1}X_t)_{t=1,\ldots,r_n}$, $\Psi_f(N)$ denotes the Laplace
functional of the point process $N$ evaluated at the non-negative
function $f$ and
$(a_n)$ satisfies $P(|X|>a_n)\sim n^{-1}$. Davis and Hsing
\cite{davis:hsing:1995} required \eqref{eq:aofa} to hold
only for non-negative measurable step functions $f$, which have bounded support
in $\ov\bbr_0$.
The mixing condition ${\mathcal A}(a_n)$ is very general.
It ensures that $N_n$ has the same limit (provided
that it exists) as
a sum of $[n/r_n]$ i.i.d. copies of the point process $N_{n,r_n}$.
Condition~${\mathcal A}(a_n)$
is implied by many known mixing conditions, particularly strong
mixing (see \cite{davis:hsing:1995}).

\subsection{\texorpdfstring{The Davis and Hsing \cite{davis:hsing:1995} approach}{The Davis and Hsing
[10] approach}}
Davis and Hsing presented a rather general approach
to the extremes of a strictly stationary sequence $(X_t)$. We quote their
Theorem 2.7 for further reference.
\begin{theorem}\label{thm:dh} Assume that $(X_t)$ is regularly
varying with index
$\alpha>0$
and normalization $(a_n)$ such that $P(|X|>a_n)\sim n^{-1}$, the
mixing condition ${\mathcal A}(a_n)$
is satisfied,
and the anticlustering condition
%
\begin{equation}\label{eq:ac}
\lim_{m\to\infty} \limsup_{\nto} P\Bigl(\max_{m\le|t|\le r_n} |X_t|>y
a_n\big\vert|X_0|>y a_n\Bigr)=0,\qquad y>0 ,
\end{equation}
holds. Here $(r_n)$ is an integer sequence such that $r_n\to\infty$,
$r_n=\mathrm{o}(n)$, which
appears in the definition of ${\mathcal A}(a_n)$. Then the following
point process convergence holds in $M_p(\ov\bbr_0)$, the set of point processes
with state space
$\ov\bbr_0$, equipped with the vague topology and
the Borel $\sigma$-field:
\[
N_n=\sum_{i=1}^n \vep_{X_t/a_n}\std N=\sum_{i=1}^\infty
\sum_{j=1}^\infty\vep_{P_i Q_{ij}} ,
\]
where $(P_i)$ are the points of a Poisson process on $(0,\infty)$
with intensity $\la(\mathrm{d}x)=\theta_{|X|} \alpha x^{-\alpha-1} \,\mathrm{d}x$ and
$\sum_{i=1}^\infty\vep_{Q_{ij}}$, $j\ge1$, constitute an i.i.d.
sequence of point processes whose points satisfy the property $|Q_{ij}|\le1$
a.s. and
$\sup_j |Q_{ij}|=1$ a.s. Here $\theta_{|X|}\in[ 0,1]$ is
the extremal index of the sequence $(|X_t|)$.
\end{theorem}
\begin{remark}\label{rem:1}
The anticlustering condition
\eqref{eq:ac} ensures that clusters of extremes become separated
from one another through time. (For a precise description of the
distribution of the point processes $\sum_{i=1}^\infty\vep_{Q_{ij}}$, see
\cite{davis:hsing:1995}. For more on the extremal index of a
stationary sequence, see \cite{leadbetter:lindgren:rootzen:1983} and
\cite{embrechts:kluppelberg:mikosch:1997}, Section 8.1. For an
introduction to point processes and their convergence in the context
of extreme value theory,
see \cite{resnick:1987,resnick:2007}.)
\end{remark}

An immediate consequence of Theorem~\ref{thm:dh} is limit theory for
the maxima $M_{n}^{|X|}=\max_{t=1,\ldots,n}|X_t|$, $n\ge1$, of the
sequence $(|X_t|)$.
Indeed, we conclude with $(a_n)$ chosen such that $nP(|X|>a_n)\sim1$,
\begin{eqnarray*}
\lim_{\nto}P\bigl(a_n^{-1} M_{n}^{|X|}\le x\bigr)&=&\lim_{\nto}P\bigl(N_n([-x,x])=0\bigr)=
P\bigl(N([-x,x])=0\bigr)
\\
&=& P\Bigl(\sup_{i\ge1}P_i\sup_{j\ge1} |Q_{ij}|\le x\Bigr)
\\
&=&P\Bigl(\sup_{i\ge
1}P_i\le x\Bigr)=P(P_1\le x)
\\
&=&\Phi_\alpha^{\theta_{|X|}}(x) ,\qquad x>0 ,
\end{eqnarray*}
where $\Phi_\alpha(x)=\exp\{-x^{-\alpha}\}$, $x>0$, denotes the Fr\'echet
distribution function with parameter $\alpha$. Similar results
can be derived for the maxima and upper-order statistics of the
$X$-sequence, joint convergence of minima and maxima, and other results
belonging to the folklore of extreme value theory.
Theorem~\ref{thm:dh} is fundamental for an extreme value theory of
the sequence $(X_t)$,
and the results reported by \cite
{basrak:davis:mikosch:2002,basrak:segers:2009,davis:hsing:1995,davis:mikosch:1998,resnick:2007}
also show that the point process convergence can be used to derive
limit results for
sums, sample autocovariances and autocorrelations, and large deviation
results.\looseness=1
\subsection{A translation result}
Our next result states that the stochastic volatility model \eqref
{eq:sv} inherits the
properties relevant for the extremal behavior of $(X_t)$
from the volatility sequence $(\sigma_t)$.
\begin{theorem}\label{thm:main}
Consider the stochastic volatility model \eqref{eq:sv}. Assume that
$(\sigma_t)$ is
regularly varying with index $\alpha>0$, it
satisfies the strong mixing property, and
$E|Z|^{\alpha+\eps}<\infty$ for some $\eps>0$.
Then
$(X_t)$ is regularly varying with index $\alpha$ and
is strongly mixing with the same rate as $(\sigma_t)$.
If $(X_t)$ also satisfies the
anticlustering condition \eqref{eq:ac}, then
Theorem~\ref{thm:dh} applies.
\end{theorem}
\begin{pf}
The proof of the regular variation of $(X_t)$ follows from
Lemma~\ref{lem:1}. Strong mixing of $(\sigma_t)$ implies
strong mixing of $(X_t)$ with the same rate (see
page 258 in \cite{davis:mikosch:2009a}). Because we assume the
anticlustering condition
\eqref{eq:ac} for $(X_t)$, the conditions of Theorem~\ref{thm:dh} are
satisfied.
\end{pf}
\begin{remark}
If $|Z|\le M$ a.s. for some positive $M$, then the anticlustering condition
for $(X_t)$ follows trivially from that for $(\sigma_t)$. If $Z$ is
unbounded, then whether this conclusion remains
true is not obvious. However, when dealing with concrete examples, it
often is not
difficult to derive the anticlustering condition for $(X_t)$; see the
examples below.\looseness=1
\end{remark}

In what follows, we consider
three examples of regularly varying
stochastic volatility models. In all cases, the
volatility sequence $(\sigma_t)$ is stationary and dependent.
We verify the regular variation, strong mixing, and anticlustering conditions
for $(\sigma_t)$ and show that these properties are inherited by
$(X_t)$. The exponential $\operatorname{AR}(1)$ model $(\sigma_t)$ of
Section~\ref{subsec:ar} does not cause extremal clustering
of $(X_t)$ whereas a random coefficient autoregressive or
linear process structure of
$(\sigma_t)$ triggers extremal dependence in the stochastic
volatility model; see
Sections~\ref{subsec:2} and \ref{subsec:3}.

\section{Exponential $\operatorname{AR}(1)$}\label{subsec:ar}
Our first example is an exponential autoregressive process of order 1
process [we write $\operatorname{AR}(1)$] given by
%
\begin{equation}\label{eq:esv}
\sigma_t= \mathrm{e}^{Y_t} ,\qquad t\in\bbz,
\end{equation}
where $(Y_t)$ is a causal stationary $\operatorname{AR}(1)$ process
$Y_t=\varphi Y_{t-1} +\eta_t$ for some $\varphi\in(-1,1)$ and
an i.i.d. sequence $(\eta_t)$ of random variables.
\begin{example}\label{exam:egarch}
Volatility sequences of the type
\eqref{eq:esv} appear in the
EGARCH (exponential GARCH) model introduced by \cite{nelson:1991}.
In this case, $X_t=\sigma_t Z_t$, $t\in\bbz$, $(Z_t)$ is an i.i.d.
sequence and
%
\begin{equation}\label{eq:egarch}
\log(\sigma_t^2)=
\alpha_0(1-\varphi)^{-1}+\sum_{k=0}^{\infty}\varphi^k
(\gamma_0Z_{t-1-k}+\delta_0 |Z_{t-1-k}|) ,\qquad t\in\bbz
\end{equation}
for positive parameters $\alpha_0,\delta_0,\gamma_0$ and
$\varphi\in(-1,1)$.
Most often, it is assumed that $(Z_t)$ is an i.i.d. standard normal
sequence.
In that case, $\sigma$ has all moments in contrast to the situation
that we
consider in this section.
This model is close to the stochastic volatility model \eqref{eq:sv}
with an
exponential $\operatorname{AR}(1)$ volatility sequence \eqref{eq:esv}. However,
in the EGARCH model, $Z_t$ feeds into the sequence $(\sigma_{s})_{s>t}$,
and thus the $\sigma$- and $Z$-sequences are dependent.
\end{example}
%
\subsection{Mixing property}\label{subsub:mix}
It is known that $(Y_t)$, and hence $(\sigma_t)$, are $\beta$-mixing with
geometric rate if $\eta$ has a positive density in some neighborhood
of $E\eta$ (cf. \cite{doukhan:1994}, Theorem 6, page 99).
%
\subsection{Regular variation}
We introduce the following conditions:
%
\begin{eqnarray}
\label{eq:bal}P(\mathrm{e}^\eta>x)&=&x^{-\alpha}L(x) ,\qquad x>0,\\
\label{eq:bal2}P(\mathrm{e}^{\eta^-}>x)&\le& c P(\mathrm{e}^\eta>x) ,\qquad x\ge1,
\end{eqnarray}
for some $\alpha>0$, a slowly varying function $L$ and some constant $c>0$, and $\eta^-$ denotes the negative
part of $\eta$.
Here and in what follows, $c$ denotes any positive constants that are
possibly different but whose values are not of interest.
Note that these conditions are
satisfied if $\eta$ is gamma or Laplace distributed.

We first prove that $\sigma$ is regularly varying.
\begin{lemma}\label{lem:1a} Assume \eqref{eq:bal} and also
\eqref{eq:bal2} if $\varphi<0$.
Then
$E\mathrm{e}^{(\alpha+\vep)\varphi Y}<\infty$ for some small $\vep>0$, and
the following relation holds:
\begingroup
\abovedisplayskip=7pt
\belowdisplayskip=7pt
\begin{equation}\label{eq:eq:ac}
P(\sigma>x)\sim
E\mathrm{e}^{\alpha\varphi Y} P(\mathrm{e}^\eta>x)= E\mathrm{e}^{\alpha\varphi Y}
x^{-\alpha}L(x) ,\qquad \xto.
\end{equation}
\end{lemma}
\begin{pf}
Because $\sigma_t=\mathrm{e}^{\eta_t} \sigma_{t-1}^\varphi$, the
random variables
$\mathrm{e}^{\eta_t}$, $\sigma_{t-1}$ are independent and, by
\eqref{eq:bal}, $\mathrm{e}^\eta$ is
regularly varying with index $\alpha>0$, we may apply a result of
Breiman \cite{breiman:1965} to conclude that \eqref{eq:eq:ac} holds if
we can show that
there exists an $\vep>0$ such that $ E\mathrm{e}^{(\alpha+\vep) \varphi
Y}<\infty$.
We first consider the case of positive $\varphi$.
Here
\[
E\mathrm{e}^{(\alpha+\vep) \varphi Y}=
\prod_{i=1}^\infty E\mathrm{e}^{(\alpha+\vep) \varphi^i \eta} .
\]
By \eqref{eq:bal}, for every $\delta>0$, there exists an $x_0>1$
such that $P(\mathrm{e}^{\eta}>x)\le x^{-\alpha+\delta}$ for $x\ge x_0$ (so-called
\emph{Potter bounds}; see Bingham \textit{et al.} \cite
{bingham:goldie:teugels:1987}, page 25). Thus for small $\vep
,\delta>0$ such that $((\alpha-\delta)/[(\alpha+\vep)\varphi^1]-1)>0$,
\begin{eqnarray*}
E\mathrm{e}^{(\alpha+\vep) \varphi^i \eta}&\le&
x_0^{(\alpha+\vep) \varphi^i}+
\int_{x_0^{(\alpha+\vep)\varphi^i}}^\infty
P\bigl(\mathrm{e}^{(\alpha+\vep) \varphi^i \eta}>y\bigr) \,\mathrm{d}y
\\[-2pt]
&\le& x_0^{(\alpha+\vep) \varphi^i}+\bigl((\alpha-\delta)/[(\alpha
+\vep)\varphi^i]-1\bigr)^{-1} x_0^{-\alpha+\delta+(\alpha+\vep)\varphi^i}.
\end{eqnarray*}
We conclude that for small $\vep,\delta>0$, some constants $c>0$,
\begin{eqnarray*}
\prod_{i=1}^\infty E\mathrm{e}^{(\alpha+\vep) \varphi^i \eta}&\le&
\exp\Biggl\{ \sum_{i=1}^\infty\bigl[x_0^{(\alpha+\vep)
\varphi^i}-1+\bigl((\alpha-\delta)/[(\alpha+\vep)\varphi^i]-1\bigr)^{-1}
x_0^{-\alpha+\delta+(\alpha+\vep)\varphi^i} \bigr] \Biggr\}
\\[-2pt]
&\le& c \exp\Biggl\{ c \sum_{i=1}^\infty\varphi^i \Biggr\}<\infty.
\end{eqnarray*}
We next consider the case of negative $\varphi$. We observe that
\begin{eqnarray*}
E\mathrm{e}^{(\alpha+\vep) \varphi Y}\le\prod_{i=1}^\infty E\mathrm{e}^{(\alpha
+\vep) \varphi^{2i} \eta}
\prod_{i=1}^\infty E\mathrm{e}^{(\alpha+\vep) |\varphi|^{2i-1} \eta^-} .
\end{eqnarray*}
Similar calculations as before, where we exploit \eqref{eq:bal} and
\eqref{eq:bal2}, show that the right-hand side is finite for small
$\vep$.
\end{pf}
\begin{lemma}\label{lem:2} Assume the conditions of Lemma~\ref{lem:1a}.
Then the sequence $(\sigma_t)$ is regularly varying with index
$\alpha$. The limit
measure of the vector ${\boldsymbol\Sigma}_d=(\sigma_1,\ldots,\sigma
_d)'$ is
given by the following limiting relation on the Borel $\sigma$-field
of $\ov\bbr{}^d_0$:
%
\begin{equation}\label{eq:limit}
\frac{
P(x^{-1}{\boldsymbol\Sigma}_d\in\cdot)}{P(\sigma>x)}\stv\alpha\sum_{i=1}^d
\int_0^\infty y^{-\alpha-1} I_{\{ y \bfe_i\in\cdot\}} \,\mathrm{d}y , \qquad\xto,
\end{equation}
where $\bfe_i$ is the $i$th unit vector in $\bbr^d$.
\end{lemma}\eject
\endgroup
\begin{pf} We give the proof only for positive $\varphi$ and
$\eta$. Proofs
for the other cases are similar.

We observe that
\[
{\boldsymbol\Sigma}_d=\operatorname{diag}(\mathrm{e}^{\varphi Y_0},\mathrm{e}^{\varphi^2
Y_0},\ldots,\mathrm{e}^{\varphi^dY_0}) \pmatrix{
\mathrm{e}^{\eta_1}
\cr
\mathrm{e}^{\eta_2+\varphi\eta_1}
\cr
\vdots
\cr
\mathrm{e}^{\eta_d+\varphi\eta_{d-1}+\cdots+ \varphi^{d-1} \eta_1}}
=\bfA\bfB.
\]
Because $E\|\bfA\|^{\alpha+\vep}<\infty$ for small positive $\vep$ and
$\bfA$ and $\bfB$ are independent regular variation of ${\boldsymbol
\Sigma
}_d$ will follow
from Breiman's multivariate result \cite{basrak:davis:mikosch:2002} if
it can
be shown that $\bfB$ is regularly varying with index $\alpha$.
Indeed, we
will show that $\bfB$ has the
same limit measure as
\[
(\mathrm{e}^{\eta_1},\mathrm{e}^{\eta_2}
E\mathrm{e}^{\alpha\varphi\eta},
\ldots,\mathrm{e}^{\eta_d} E\mathrm{e}^{\alpha\varphi\eta}\cdots
E\mathrm{e}^{\alpha\varphi^{d-1}\eta} )'.
\]
This fact does
not follow from the continuous mapping theorem for regularly varying vectors
(see \cite{hult:lindskog:2005,hult:lindskog:2006a}), because the function
 $(r_1,\ldots,r_d)\to(r_1,r_1^\varphi r_2,\ldots,
r_1^{\varphi^{d-1}}\cdots r_{d-1}^{\varphi}r_d)$
does not have the homogeneity property.

For simplicity, we prove the result only for $d=2$, the general case
being analogous. To ease notation, we also write
$R_i=\mathrm{e}^{\eta_i}$, $i=1,2$. Choose $a_n$ such that $P(\mathrm{e}^\eta>a_n)\sim
n^{-1}$ and
take any set $A\subset\ov\bbr{}^2_0$ that is a subset of the first
orthant bounded away from 0 and
continuous with respect to  the limiting measure of ${\boldsymbol\Sigma
}_d$ in the
formulation of
the lemma. Write $B=\{a_n^{-1} (R_1,R_1^\varphi R_2)\in A\}$, and for
any $\vep,\gamma>0$, consider the disjoint sets
\begin{eqnarray*}
B_1&=&B\cap\{R_1>\vep a_n,R_2>\gamma a_n\} ,
\\
B_2&=&B\cap\{R_1>\vep a_n,R_2\le\gamma a_n\} ,
\\
B_3&=&B\cap\{R_1\le\vep a_n,R_2>\gamma a_n\} ,
\\
B_4&=&B\cap\{R_1\le\vep a_n,R_2\le\gamma a_n\} .
\end{eqnarray*}
Then for any $\vep,\gamma>0$,
\[
nP(B_1)\le n P(R_1>\vep a_n) P(R_2>\gamma a_n)\to0 .
\]
Next, consider $B_3$. Choose some $M>1$ and consider the disjoint
partition
of $B_3$,
\[
B_{31}=B_3\cap\{R_1\in[1,M]\} ,\qquad
B_{32}=B_3\cap\{R_1>M\} .
\]
Then
\[
nP(B_{32})\le n P(R_2>\gamma a_n) P(R_1>M)\sim\gamma^{-\alpha}
P(R_1>M),\qquad \nto.
\]
Thus, for any $\vep,\gamma>0$,
\[
\lim_{M\to\infty}\limsup_{\nto}n P(B_{32})=0 .
\]
%
Observe that $nP(R_1I_{\{R_1\in[1,M]\}}> c a_n)\to0$ for every $c>0$
and, by Breiman's result \cite{breiman:1965}, $R_2R_1^\varphi
I_{\{R_1\in[1,M]\}}$ is regularly varying. By Lemma 3.12 of
\cite{jessen:mikosch:2006},
\[
\bigl(R_1 I_{\{R_1\in[1,M]\}}, R_2R_1^\varphi
I_{\{R_1\in[1,M]\}}\bigr)=\bigl(R_1
I_{\{R_1\in[1,M]\}}, 0\bigr)+
\bigl(0,R_2R_1^\varphi
I_{\{R_1\in[1,M]\}}\bigr)
\]
is regularly varying with the same index and limiting measure as
$(0,R_2R_1^\varphi
I_{\{R_1\in[1,M]\}})$. Therefore,
\begin{eqnarray*}
nP(B_{31})&\sim& nP \bigl(a_n^{-1}\bigl(0, R_2R_1^\varphi
I_{\{R_1\in[1,M]\}}\bigr)\in A , R_2>\gamma a_n \bigr)
\\
&=&nP \bigl(a_n^{-1}R_2 R_1^\varphi
I_{\{R_1\in[1,M]\}} \in\operatorname{proj}_2 A , R_2>\gamma a_n
\bigr) I_{\{\operatorname{proj}_1A=\{0\}\}} ,
\end{eqnarray*}
where $\operatorname{proj}_i A$, $i=1,2$, are the projections of $A$ on the
$x$- and $y$-axes, respectively. Regular variation of $R_2$ with limit
measure $\mu(t,\infty)=
t^{-\alpha}$, $t>0$,
ensures that
\begin{eqnarray*}
\lim_{M\to\infty}\lim_{\gamma\to0}\lim_{\nto} n P(B_{31})&=&
\lim_{M\to\infty}\lim_{\gamma\to0}
E\mu\bigl\{t>\gamma\dvt R_1^\varphi I_{\{R_1\in[1,M]\}} t \in\operatorname{proj}_2
A\bigr\} I_{\{\operatorname{proj}_1A=\{0\}\}}
\\
&=&\lim_{M\to\infty}
E\mu\bigl\{t>0\dvt R_1^\varphi I_{\{R_1\in[1,M]\}} t \in\operatorname{proj}_2
A\bigr\} I_{\{\operatorname{proj}_1A=\{0\}\}}
\\
&=& \lim_{M\to\infty} ER_1^{\alpha\varphi} I_{\{R_1\in[1,M]\}}
\mu(\operatorname{proj}_2 A) I_{\{\operatorname{proj}_1A=\{0\}\}}
\\
&=& ER_1^{\alpha\varphi} \mu(\operatorname{proj}_2 A) I_{\{
\operatorname{proj}_1A=\{0\}\}} .
\end{eqnarray*}
We have $A\subset\{{\mathbf x}\dvt|x_1|+|x_2|>\delta\}$ for small $\delta>0$.
Then
$B_4$ is contained in the
union of the following sets for $M>1$:
\begin{eqnarray*}
B_{41}&=&B_4\cap\{R_1>0.5 \delta a_n\} ,
\\
B_{42}&=&B_4\cap\{R_1^\varphi R_2>0.5 \delta a_n,R_1> M\} ,
\\
B_{43}&=&B_4\cap\{R_1^\varphi R_2>0.5 \delta a_n,R_1\in[1,M]\} .
\end{eqnarray*}
Choosing $\vep$ sufficiently small, $B_{41}$ is empty.
Moreover,
by Breiman's result,
\[
nP(B_{42})\le n P(R_1^\varphi R_2>0.5 \delta a_n,R_1> M)\sim c
E\bigl[R_1^{\alpha\varphi}I_{\{R_1>M\}}\bigr] ,\qquad \nto.
\]
Choosing $\gamma$ sufficiently small, the set $B_{43}$ is empty.
Therefore, and because\break $E[R_1^{(\alpha+\vep)\varphi}]<\infty$,
\[
\lim_{M\to\infty}\limsup_{\nto} n P(B_{4i})=0 , \qquad i=1,2,3.
\]
It remains to consider the set $B_2$. Consider the disjoint partition
of $B_2$ for $M>1$,
\[
B_{21}=B_2\cap\{R_2\le M\} \quad\mbox{and}\quad B_{22}=B_2\cap\{R_2>M\} .
\]
Because $P(B_{22})\le P(R_1>\vep a_n)P(R_2>M)$, we have
\[
\lim_{M\to\infty}\limsup_{\nto} n P(B_{22})=0 .
\]
Moreover,
\[
n P(B_{21})\sim n P\bigl(R_1>\vep a_n , a_n^{-1}(R_1,0)\in A\bigr) .
\]
Thus, for every $M>0$,
\begin{eqnarray*}
\lim_{\vep\to0}
\lim_{\nto} n P(B_{21})&=&
\lim_{\vep\to0} \mu\{t>\vep\dvt t\in\operatorname{proj}_1 A\} I_{\{
\operatorname{proj}_2 A=\{0\} \}}
\\
&=& \mu\{\operatorname{proj}_1 A\} I_{\{ \operatorname{proj}_2A=\{0\} \}} .
\end{eqnarray*}
Summarizing the foregoing arguments, we have proven that
\[
n P\bigl(a_n^{-1}(R_1,R_1^\varphi R_2)\in A\bigr)\to\alpha\int_0^\infty
x^{-\alpha-1} \bigl[I_{\{x\bfe_1 \in A\}}+ ER_1^{\alpha\varphi} I_{\{
x\bfe_2\in A\}}\bigr] \,\mathrm{d}x .
\]
Modifying the proof above for $d\ge2$, we obtain
\begin{eqnarray*}
&&n P(a_n^{-1}\bfB\in A)
\\
&&\quad\to \alpha\int_0^\infty
x^{-\alpha-1} \bigl[I_{\{x\bfe_1 \in A\}}+
E\mathrm{e}^{\alpha\varphi\eta} I_{\{x\bfe_2\in A\}}+\cdots
+E\mathrm{e}^{\alpha\varphi\eta} \cdots E\mathrm{e}^{\alpha\varphi^{d-1}\eta}
I_{\{x\bfe_d\in A\}} \bigr] \,\mathrm{d}x .
\end{eqnarray*}
We now apply the multivariate Breiman result
\cite{basrak:davis:mikosch:2002} to obtain
\begin{eqnarray*}
&&n P(a_n^{-1}\bfA\bfB\in A)
\\
&&\quad\to \alpha\int_0^\infty
x^{-\alpha-1} E \bigl[I_{\{ \mathrm{e}^{\varphi Y}x\bfe_1 \in A\}}+
E\mathrm{e}^{\alpha\varphi\eta} I_{\{\mathrm{e}^{\varphi^2 Y }x\bfe_2\in A\}}+\cdots
\\
&&\hphantom{\quad\to \alpha\int_0^\infty
x^{-\alpha-1} E \bigl[}{}+E\mathrm{e}^{\alpha\varphi\eta} \cdots E\mathrm{e}^{\alpha\varphi^{d-1}\eta}
I_{\{\mathrm{e}^{\varphi^d Y}x\bfe_d\in A\}} \bigr] \,\mathrm{d}x
\\
&&\quad=
\alpha\int_0^\infty
x^{-\alpha-1} \bigl[ E\mathrm{e}^{\alpha\varphi Y} I_{\{x\bfe_1 \in A\}}+
E\mathrm{e}^{\alpha\varphi\eta}E\mathrm{e}^{\alpha\varphi^2Y} I_{\{x\bfe_2\in A\}
}+\cdots
\\
&&\hphantom{\quad=
\alpha\int_0^\infty
x^{-\alpha-1} \bigl[}{}+E\mathrm{e}^{\alpha\varphi\eta} \cdots E\mathrm{e}^{\alpha\varphi^{d-1}\eta}
E\mathrm{e}^{\alpha\varphi^{d}Y}
I_{\{x\bfe_d\in A\}} \bigr] \,\mathrm{d}x
\\
&&\quad=\alpha E\mathrm{e}^{\alpha\varphi Y}\int_0^\infty
x^{-\alpha-1} \bigl[ I_{\{x\bfe_1 \in A\}}+
I_{\{x\bfe_2\in A\}}+\cdots
+
I_{\{x\bfe_d\in A\}}\bigr] \,\mathrm{d}x .
\end{eqnarray*}
This relation and Lemma~\ref{lem:1a} conclude the proof for $\varphi
\in
(0,1)$ and $\eta>0$ a.s.
\end{pf}
%
\subsection{Anticlustering condition}
\begin{lemma}\label{lem:acsv} Assume \eqref{eq:bal} and also
%
\begin{equation}\label{eq:bal3}
P\bigl(\mathrm{e}^{|\eta|}>x\bigr)\le
c P(\mathrm{e}^\eta>x) ,\qquad x\ge1
\end{equation}
for some $c>0$, $\varphi\in(-1,1)$. Then the
anticlustering condition \eqref{eq:ac} holds for the sequence $(\sigma_t)$ and any sequence $(r_n)$ satisfying $r_n=\mathrm{O}(n^{\gamma})$
for some
$\gamma\in(0,1)$.
If $|Z|$ has finite moments of any order,
then \eqref{eq:ac} is also satisfied for the stochastic volatility sequence $(X_t)$ with
the same sequence $(r_n)$ as for $(\sigma_t)$.
If $E|Z|^{\alpha+\xi}<\infty$ for some $\xi>0$,
then \eqref{eq:ac} holds for the sequence $(X_t)$ with $(r_n)$ such
that $r_n = \mathrm{O}(n^\gamma)$ for every $\gamma\in(0,1)$.
\end{lemma}
\begin{pf} Throughout, we assume that $\varphi\ne0$. If
$\varphi=0$,
then both $(\sigma_t)$ and $(X_t)$ are
i.i.d. regularly varying sequences, and \eqref{eq:ac} is trivially satisfied.

We first prove the result for $(\sigma_t)$. We begin under the
assumptions $\varphi\in(0,1)$ and $\eta>0$, and verify that
%
\begin{equation}\label{eq:2}
\lim_{m\to\infty}\limsup_{\nto}
P\Bigl(\max_{m\le t\le r_n} Y_t>\log(y a_n) \big\vert Y_0>\log
(ya_n)\Bigr)=0 , \qquad y>0 .
\end{equation}
Fix $y>0$ and write
$
B=\{\max_{m\le t\le r_n} Y_t>\log(ya_n) \}$ and observe that
%
\begin{equation}\label{eq:3}
Y_t=\varphi^{t-m}Y_m +\sum_{i=m+1}^t\varphi^{t-i} \eta_i , \qquad m\le t .
\end{equation}
Then $B\subset B_1\cup B_2$, where for $\delta\in(0,1)$,
\[
B_1=\{Y_m>\delta\log(ya_n)\} \quad\mbox{and}\quad
B_2= \Biggl\{\max_{m\le t\le r_n} \sum_{i=m+1}^t\varphi^{t-i} \eta_i
>(1-\delta)\log(ya_n) \Biggr\} .
\]
Because $Y_0$ is independent of $(\eta_t)_{t\ge1}$,
$P(B_2)=P(B_2\vert Y_0>\log(y a_n))$. Therefore, and
by Markov's inequality,
\begin{eqnarray*}
P\bigl(B_2\vert Y_0>\log(y a_n)\bigr)
&\le&\sum_{t=m}^{r_n}P \Biggl( \sum_{i=m+1}^t\varphi^{t-i} \eta_i
>(1-\delta)\log(y a_n) \Biggr)
\\
 &\le&\sum_{t=m}^{r_n}P\bigl( Y_t
>(1-\delta)\log(y a_n)\bigr)
\\
&\le&r_n P\bigl(Y>(1-\delta)\log(ya_n)\bigr)
\\
&\le& r_n E\mathrm{e}^{(\alpha-\vep)Y} (ya_n)^{-(1-\delta) (\alpha
-\vep)}
\end{eqnarray*}
for $0<\vep<\alpha$ and large $n$.
Because $E\mathrm{e}^{(\alpha-\vep)Y} <\infty$ and $a_n=n^{1/\alpha
}\ell(n)$
for some slowly varying function $\ell$, choosing $\delta,\vep>0$
sufficiently
small, the right-hand side converges to 0 if
$r_n=\mathrm{O}(n^{\gamma})$ for some $\gamma<1$. Moreover, $B_1\subset
B_{11}\cup B_{12}$, where
\[
B_{11}=\{\varphi^m Y_0>0.5\delta\log(ya_n)\} \quad\mbox{and}\quad
B_{12}=\Biggl\{\sum_{i=1}^m\varphi^{m-i} \eta_i> 0.5\delta\log(ya_n)\Biggr\} .
\]
For any $m$, small $\vep>0$, large $n$, we have
\[
nP\bigl(B_{11}\cap\{Y_0>\log(ya_n)\}\bigr)=nP\bigl( Y_0>0.5\delta\log(ya_n)
\varphi^{-m}\bigr)
\le n (ya_n)^{-0.5 \delta\varphi^{-m}(\alpha-\vep)} .
\]
Therefore, choosing $m$ sufficiently large, the right-hand side converges
to 0.
Because $Y_0$ and $B_{12}$ are independent,
\[
P\bigl(B_{12}\vert Y_0>\log(ya_n)\bigr)= P(B_{12}) .
\]
The right-hand side is bounded by $P(Y> 0.5\delta\log(ya_n))=\mathrm{o}(1)$.
Thus, we have proven
\[
\lim_{\nto} P\bigl(B_1\vert Y_0>\log(y a_n)\bigr)=0
\]
and that \eqref{eq:2} holds. Next, we prove
%
\begin{equation}\label{eq:4}
\lim_{m\to\infty}\limsup_{\nto}
P\Bigl(\max_{-r_n\le t\le-m} Y_t>\log(ya_n) \big\vert Y_0>\log(y a_n)\Bigr)=0 .
\end{equation}
Write
\[
C=\Bigl\{\max_{-r_n\le t\le-m} Y_t>\log(ya_n) , Y_0>\log(ya_n)\Bigr\} .
\]
Again using \eqref{eq:3}, we see that
$C\subset C_1\cup C_2$,
where, for $\delta\in(0,1)$,
\begin{eqnarray*}
C_1&=&\{Y_{-r_n}>\delta\log a_n,Y_0>\log(ya_n) \} ,
\\
C_2&=& \Biggl\{\max_{-r_n\le t\le-m }
\sum_{i=-r_n+1}^t\varphi^{t-i}\eta_i>(1-\delta) \log(ya_n)
,Y_0>\log(ya_n) \Biggr\}.
\end{eqnarray*}
Another application of \eqref{eq:3} and stationarity yields
\begin{eqnarray*}
n P(C_1)&\le& n P\bigl(Y_0>\delta\log(ya_n) ,Y_0> (1-\delta)
\varphi^{-r_n}\log(ya_n)\bigr)
\\
&&{}+n P \Biggl(Y_0>\delta\log(ya_n) ,\sum_{i=1}^{r_n}
\varphi^{r_n-i}\eta_i>\delta
\log(ya_n) \Biggr)=I_1+I_2 .
\end{eqnarray*}
By regular variation, for small $0<\vep<\alpha$ and large $n$,
\[
I_1\le n (ya_n)^{-(\alpha-\vep) (1-\delta) \varphi^{-r_n}} .
\]
Because $r_n\to\infty$, we have $I_1=\mathrm{o}(1)$ as $\nto$. Moreover, it
follows that
\[
\limsup_{\nto}I_2\le c \limsup_{\nto}P \Biggl(\sum_{i=1}^{r_n}
\varphi^{r_n-i}\eta_i>\delta
\log(ya_n)
\Biggr)=0 .
\]
Thus, we have proven that $\limsup_{\nto} n P(C_1)=0$. For $C_2$, we have
$C_2\subset C_{21}\cup C_{22}$, where
\begin{eqnarray*}
C_{21}&=& \Biggl\{\max_{-r_n\le t\le-m }
\sum_{i=-r_n+1}^t\varphi^{t-i}\eta_i>(1-\delta) \log(ya_n)
,\varphi^m
Y_{-m}>\delta\log(ya_n) \Biggr\} ,
\\
C_{22}&=& \Biggl\{\max_{-r_n\le t\le-m }
\sum_{i=-r_n+1}^t\varphi^{t-i}\eta_i>(1-\delta) \log(ya_n) ,\sum
_{i=-m+1}^0\varphi^{-i}\eta_i>
(1-\delta) \log(ya_n) \Biggr\} .
\end{eqnarray*}
Thus, for small $\vep$, large $m$,
\[
n P(C_{21})\le n P\bigl(Y_0>\varphi^{-m} \delta\log(y a_n)\bigr)\le n
(ya_n)^{-(\alpha-\vep)\delta\varphi^{-m}}\to0 , \qquad\nto,
\]
and for small $\vep,\delta$,
\begin{eqnarray*}
nP(C_{22})&=&n P \Biggl(\max_{-r_n\le t\le-m }
\sum_{i=-r_n+1}^t\varphi^{t-i}\eta_i>(1-\delta) \log(ya_n) \Biggr)
\\
&&{}\times P \Biggl(\sum_{i=-m+1}^0\varphi^{-i}\eta_i>
(1-\delta) \log(ya_n) \Biggr)
\\
&\le&n r_n \bigl[P\bigl(Y>(1-\delta) \log(ya_n)\bigr)\bigr]^2\le n r_n (ya_n)^{-2
(\alpha-\vep)(1-\delta)}=\mathrm{o}(1),
\end{eqnarray*}
provided that $r_n=\mathrm{O}(n^{\gamma})$ for some $\gamma<1$.
For general $\eta$ and $|\varphi|<1$, we
see that $|Y_t|\le\sum_{i=j} |\varphi|^j |\eta_{t-j}|$. We
can apply the same reasoning as above, using \eqref{eq:bal3}.

We now turn to the proof of the anticlustering condition for
$(X_t)$. An inspection of the foregoing proof shows that we have to add
the terms $R_t=\log|Z_t|$ to $|Y_t|$. We restrict ourselves to the
cases $\varphi\in(0,1)$, $\eta>0$ a.s., and only show
that
%
\begin{equation}\label{eq:2+}
\lim_{m\to\infty}\limsup_{\nto}
P\Bigl(\max_{m\le t\le r_n} (Y_t+R_t)>\log(y a_n) \big\vert Y_0+R_0>\log
(ya_n)\Bigr)=0 , \qquad y>0 .
\end{equation}
We use the same notation for the modified events. We start by
observing that
\[
B=\Bigl\{\max_{m\le t\le r_n} (Y_t+R_t)>\log(ya_n) \Bigr\}\subset B_1\cup B_2 ,
\]
where $B_2$ is the same as above and
\begin{eqnarray*}
B_1&=&\Bigl\{Y_m+\max_{m\le t\le r_n} R_t>\delta\log(ya_n)\Bigr\}
\\
&\subset&
\{Y_m>0.5 \delta\log(ya_n)\}\cup\Bigl\{\max_{m\le t\le r_n}
R_t>0.5\delta\log(ya_n)\Bigr\}=D_1\cup D_2 .
\end{eqnarray*}
Now $P(D_1)$ can be treated in the same way as $P(B_1)$ in the
foregoing proof. If $|Z|$
has moments of any order $h>0$, then an application of Markov's
inequality for sufficiently large $h$ yields, for any choice of $r_n=o(n)$,
\[
P(D_2)=P\Bigl(\max_{m\le t\le r_n} |Z_t| >(ya_n)^{0.5\delta}\Bigr)\le r_n
P\bigl(|Z|>(ya_n)^{0.5\delta}\bigr)\le c r_n (ya_n)^{-0.5 h\delta}=\mathrm{o}(1) .
\]
On the other hand, if $r_n=\mathrm{O}(n^{\gamma})$ for every small $\gamma$,
then Markov's inequality of
order $h=\alpha+\xi$ yields the same result by choosing $\gamma$ close
to 0. This completes the proof
of \eqref{eq:2+}.
\end{pf}
%
\subsection{Main result for the exponential $\operatorname{AR}(1)$ process}
Here we give sufficient conditions for the validity of
Theorem~\ref{thm:dh} when $(X_t)$ is a stochastic volatility process
and the volatility
process $(\sigma_t)$ is an exponential $\operatorname{AR}(1)$ process.
The result is a consequence of the translation
result Theorem~\ref{thm:main} and the foregoing calculations.

\begin{theorem}\label{thm:mainsv}
Consider the stochastic volatility model \eqref{eq:sv}, where the
volatility sequence $(\sigma_t)$
is an exponential $\operatorname{AR}(1)$ process \eqref{eq:esv} for some $\varphi\in
(-1,1)$. Assume the following conditions:
\begin{itemize}
\item
The regular variation conditions \eqref{eq:bal} and \eqref{eq:bal3} hold for some index
$\alpha>0$.
\item
The random variable $\eta$ has positive density in some neighborhood
of $E\eta$.
\end{itemize}
Then the following properties hold for $(\sigma_t)$.
\begin{enumerate}[(3)]
\item[(1)]
Regular variation with index $\alpha$ and limiting measures
given in \eqref{eq:limit}.
\item[(2)]
$\beta$-mixing with
geometric rate and ${\mathcal A}(a_n)$ are satisfied for any sequence $(r_n)$ satisfying $r_n\ge c\log n$ for some $c>0$ and $r_n=\mathrm{o}(n)$.
\item[(3)]
The
anticlustering condition for
$r_n=\mathrm{O}(n^{\gamma})$ for any $\gamma\in(0,1)$.
\end{enumerate}
The following properties hold for the stochastic volatility process $(X_t)$:
\begin{enumerate}[(4)]
\item[(4)]
The strong mixing property with geometric
rate and ${\mathcal A}(a_n)$ are satisfied for any sequence $(r_n)$ satisfying $r_n\ge c\log n$ for some $c>0$ and
$r_n=\mathrm{o}(n)$.
\end{enumerate}
Also assume that
\begin{itemize}
\item
$E|Z|^{\alpha+\delta}<\infty$ for some $\delta>0$.
\end{itemize}
Then
\begin{enumerate}[(6)]
\item[(5)]
$(X_t)$ is regularly varying with index $\alpha$ and limiting
measures given
in Example~\ref{exam:limit}.
\item[(6)]
$(X_t)$ satisfies the anticlustering condition \eqref{eq:ac} for
$(r_n)$ such that $r_n=\mathrm{O}(n^\gamma)$ for every $\gamma<1$.
\end{enumerate}
Moreover, if
\begin{itemize}
\item
$Z$ has all moments,
\end{itemize}
then
\begin{enumerate}[(7)]
\item[(7)] $(X_t)$ satisfies the anticlustering condition
\eqref{eq:ac} for any sequence $(r_n)$ such that $r_n=\mathrm{O}(n^\gamma)$
for some
$\gamma<1$.
\end{enumerate}
In particular, Theorem~\ref{thm:dh} applies to the sequences $(\sigma_t)$
and $(X_t)$.
\end{theorem}
\begin{pf} We first give the proof for the volatility sequence $(\sigma_t)$.
Regular variation of $(\sigma_t)$ follows from Lemma~\ref{lem:2}, and
$\beta$-mixing with geometric rate follows from
Section~\ref{subsub:mix}. It follows from
\cite{davis:hsing:1995} and references therein that condition
${\mathcal A}(a_n)$ is satisfied with $r_n\ge c\log n$ for some $c>0$.
Condition \eqref{eq:ac} for $(\sigma_t)$ follows from
Lemma~\ref{lem:acsv} under the assumption that $r_n=\mathrm{O}(n^{\gamma})$
for some $\gamma\in(0,1)$.\eject

Because $\beta$-mixing with geometric rate implies strong mixing with
geometric rate and, using the argument on page 258 of
\cite{davis:mikosch:2009a}, it
follows that $(X_t)$ is strongly mixing with geometric rate. It follows from
\cite{davis:hsing:1995} and references therein that condition
${\mathcal A}(a_n)$ is satisfied for any $r_n\ge c\log n$ for some
$c>0$.
Regular variation of $(X_t)$ follows from Theorem~\ref{thm:main},
and the limiting measures are derived in
Example~\ref{exam:limit}. Finally, condition \eqref{eq:ac} was
verified in Lemma~\ref{lem:acsv}.
\end{pf}

\begin{figure}

\includegraphics{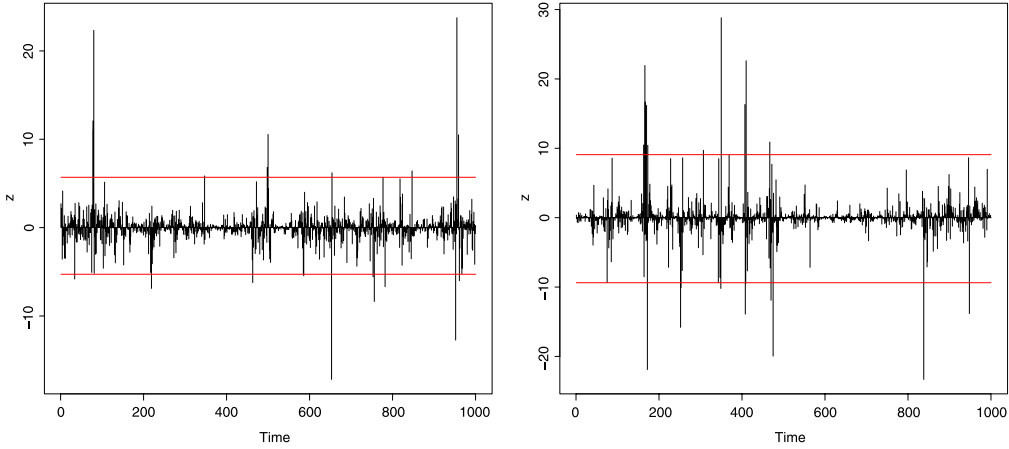}

  \caption{1000 realizations of a stochastic volatility model,
where $(\log
\sigma_t)$
is an $\operatorname{AR}(1)$ process with $\varphi=0.9$. The parallel lines indicate the
0.01 and 0.99 quantiles of the distribution of $X$. \emph{Left:}
The random variable $\eta$ is Laplace distributed: $P(X>x)=P(X
\le-x)=0.5\mathrm{e}^{-4x}$, $x>0$, and $Z$ standard normal. \emph{Right:}
The random variable $\eta$ is
$N(0,0.25)$-distributed and
$Z$ is $t$-distributed with 4 degrees of freedom
standardized to unit variance. In both graphs, $(X_t)$ is regularly
varying with index 4, and there is no extremal clustering in the sense that
high and low exceedances of the lines occur separated through time.}\label{fig:1}
\end{figure}

Using the machinery in
\cite{basrak:davis:mikosch:2002,davis:hsing:1995,davis:mikosch:1998,davis:mikosch:2001},
we can now derive various limit
results for the sequence $(X_t)$. These include infinite variance limits
for the normalized partial sums $\sum_{t=1}^n X_t$ and sample covariances
$\sum_{t=1}^{n-h} X_{t}X_{t+h}$ in the case where $\alpha<2$. For general
$\alpha>0$, the fact that the limit measures of the regular variation of $(X_t)$ are
concentrated on the axes implies that the normalized partial maxima
of $(X_t)$ converge to a Fr\'echet distribution,
%
\begin{equation}\label{eq:nocl}
\lim_{\nto}P\Bigl(a_n^{-1} \max_{t=1,\ldots,n} X_t\le
x\Bigr)=\Phi_\alpha(x)=\mathrm{e}^{-p x^{-\alpha}} ,\qquad x>0 ,
\end{equation}
where $(a_n)$ satisfies $n P(|X|>a_n)\to1$ and
\[
\lim_{\xto}\frac{P(X>x)}{P(|X|>x)}=\frac{EZ_+^\alpha}{E|Z|^\alpha
}=p\in[0,1] .
\]
Relation \eqref{eq:nocl} means that the extremal index of the
sequence $(X_t)$ is 1; that is, we get the same result as for an i.i.d.
sequence $(\wt
X_t)$ with $\wt X\eqd X$. In other words, the stochastic volatility model does not
exhibit extremal clustering. This is analogous to stochastic
volatility models
in which $E\sigma^{\alpha+\delta}<\infty$ and $Z$ is regularly
varying with index
$\alpha$ (see
\cite{davis:mikosch:2001,davis:mikosch:2009}),
although the reasons are very different in the two cases.
Figure~\ref{fig:1} presents graphs of regularly varying stochastic
volatility models with
light-tailed and heavy-tailed multiplicative noise.
In the
present case,
the structure of the
limiting measures for the regularly varying finite-dimensional
distributions of the $\sigma$-sequence is responsible for the limiting measures of the $X$-sequence.

In passing, we mention that
a condition of type \eqref{eq:bal} limits the choice of the
distributions of the
noise variable $\eta$ in the
exponential $\operatorname{AR}(1)$ process. If $\eta$ has a slightly heavier right tail than
suggested by \eqref{eq:bal} the random variable $Y$ will not have any
moments. This occurs,
for example, when $\eta$ has a lognormal or Student distribution. Thus
regular variation
 of $(\sigma_t)$ and $(X_t)$ is
possible only for a relatively thin class of noise variables $\eta$.

Before we consider other stochastic volatility models with genuine
extremal clustering, we show that the EGARCH model from
Example~\ref{exam:egarch} is regularly varying and does not have extremal
clusters.
\begin{example}
Recall the definition of the EGARCH model from
Example~\ref{exam:egarch}, particularly the dynamics of
$(\sigma_t^2)$ given by \eqref{eq:egarch}. Writing
$\eta_t=0.5 (\alpha_0(1-\varphi)^{-1}+\gamma_0Z_{t}+\delta_0
|Z_t|)$ and assuming the
conditions of Lemma~\ref{lem:1a}, we conclude that $(\sigma_t)$ is
regularly varying with index $\alpha$, and the limiting measures are
concentrated
on the axes. Using the modified Breiman lemma
from \cite{jessen:mikosch:2006}, an inspection of the proof of
Lemma~\ref{lem:2} shows that ${\boldsymbol\Sigma}_d=(\sigma_1,\ldots
,\sigma_d)'$
and $(\mathrm{e}^{\eta_0},\ldots,\mathrm{e}^{\eta_{d-1}})'E\mathrm{e}^{\alpha\varphi\eta}$ have
the same limit measures of regular variation. Therefore, regular
variation of
$\bfX_d=(X_1,\ldots,X_d)'$ will follow if we can show that
$\bfR_d=(Z_1 \mathrm{e}^{\eta_0},\ldots,Z_d \mathrm{e}^{\eta_{d-1}})'$ is regularly
varying with limit
measures concentrated on the axes. By Breiman's result, $Z_1 \mathrm{e}^{\eta_0}$
is regularly varying with index $\alpha$. Let $(a_n)$ be such that $n P(\mathrm{e}^\eta>a_n)\to1$. By construction, $Z$ has all moments, and
thus we can
choose a sequence $c_n\to\infty$
such that $n P(|Z|>c_n)\to0$ and $a_n/c_n\to\infty$. Then, for
$d\ge2$,
$\delta>0$,
\begin{eqnarray*}
n P(|Z_i \mathrm{e}^{\eta_{i-1}}|>\delta a_n,i=1,\ldots,d)&\le&
n P(|Z_1 \mathrm{e}^{\eta_{0}}|>\delta a_n,|Z_2 \mathrm{e}^{\eta_{1}}|>\delta
a_n)
\\
&\le&n P(|Z_1 \mathrm{e}^{\eta_{0}}|>\delta a_n) P(\mathrm{e}^{\eta_{1}}>\delta
a_n/c_n)+ n P(|Z|>c_n)
\\
&=&\mathrm{o}(1) .
\end{eqnarray*}
Thus, if $n P(a_n^{-1}\bfR_d \in\cdot)$ has a non-vanishing vague
limit, then it must be
concentrated on the axes. To show this, we focus on the case where
$d=2$. Here,
for $x,\delta>0$, by Breiman's result and the previous calculations,
\begin{eqnarray*}
n P(a_n^{-1} |Z_1| \mathrm{e}^{\eta_0}\le\delta,a_n^{-1} Z_2 \mathrm{e}^{\eta_1}>
x)
&=& n P(a_n^{-1}Z_2 \mathrm{e}^{\eta_1}>x)
\\
&&{}- n P(a_n^{-1}|Z_1| \mathrm{e}^{\eta_0}>
\delta, a_n^{-1}Z_2 \mathrm{e}^{\eta_1}>x)
\\
&\sim& x^{-\alpha} EZ_+^{\alpha} ,
\\
nP(a_n^{-1} |Z_1|\mathrm{e}^{\eta_0}\le\delta,a_n^{-1} Z_2\mathrm{e}^{\eta_1}\le-
x)
&=& n P(a_n^{-1}Z_2 \mathrm{e}^{\eta_1}\le-x)
\\
&&{}- nP(a_n^{-1}|Z_1|\mathrm{e}^{\eta_0}>
\delta, a_n^{-1}Z_2\mathrm{e}^{\eta_1}\le-x)
\\
&\sim& x^{-\alpha} EZ_-^{\alpha} .
\end{eqnarray*}
Therefore, $(X_t)$ is regularly varying with index $\alpha$, and the
limiting measures
are concentrated on the axes. Furthermore, if $\eta$ has a positive
density in some neighborhood of $E\eta$, then $(\log\sigma_t)$, hence
$(X_t)$, is strongly mixing with geometric rate, and then ${\mathcal A}(a_n)$
holds for any sequence $(r_n)$ satisfying $r_n=\mathrm{o}(n)$ and $r_n\ge c\log n$
for some
$c>0$. For the proof of the anticlustering condition of $(X_t)$, we
can follow along the lines of the proof of Lemma~\ref{eq:bal3}, observing
that $Z$ has all moments. Thus, the conditions of Theorem~\ref{thm:dh}
are satisfied, in particular because the limiting measures of the
regularly varying finite-dimensional distributions
of $(X_t)$ are concentrated on the axes the extremal index
$\theta_{|X|}=1$, that is, there is no extremal clustering in this
sequence.
\end{example}

\section{Stochastic recurrence equations}\label{subsec:2}
We assume that the stationary sequence $(\sigma_t)$
satisfies the relation
%
\begin{equation}\label{eq:5}
\sigma_t^p= A_t\sigma_{t-1}^p+ B_t , \qquad t\in\bbz
\end{equation}
for an i.i.d. sequence $((A_t,B_t))_{t\in\bbz}$ of non-negative
random variables and
some positive $p$. Throughout we assume the conditions of Kesten
\cite{kesten:1973}, which ensure that \eqref{eq:5} has a strictly stationary
solution, namely $E\log A<0$ and $E\log^+ B<\infty$.
\begin{example}\label{exam:garch}
For $p=2$, a model of the type
\eqref{eq:5} has attracted major
attention in the financial time series literature
\cite{andersen:davis:kreiss:mikosch:2009}: the GARCH process of
order $(1,1)$ (we write \garch) given by
$\wt X_t=\sigma_t \eta_t$, $t\in\bbz$, $(\eta_t)$
is an i.i.d. centered sequence with unit variance and
$\sigma_t^2=\alpha_0+\sigma_{t-1}^2 (\alpha_1
\eta_{t-1}^2+\beta_1)$ for positive parameters $\alpha_i,\beta_1$.
The main difference from the stochastic volatility model \eqref
{eq:sv} with the same
sequence $(\sigma_t)$ is that $\eta_t$ feeds into $(\sigma
_{t+k})_{k\ge1}$,
and thus the noise $(\eta_t)$ and $(\sigma_t)$ are dependent.
\end{example}
%
\subsection{Mixing property}
It follows from \cite{mokkadem:1990} that
$(\sigma_t^p)$ is strongly mixing with geometric rate if $A,B$ satisfy
some regularity condition. In particular, if $A_t$ and $B_t$ are
polynomials of an i.i.d. sequence $(\eta_t)$ and $\eta$ has a positive
density in
some neighborhood of $E\eta$, then $(\sigma_t)$ is $\beta$-mixing with
geometric rate. Thus the $\operatorname{GARCH}(1,1)$ model satisfies this condition
for $p=2$ if $\eta$ has a positive density in some neighborhood of the
origin.
%
\subsection{Regular variation}
Regular variation of the marginal distribution of the solution to the
stochastic recurrence equation \eqref{eq:5} was proven by Kesten \cite{kesten:1973} and Goldie
\cite{goldie:1991}. In particular, they showed that
%
\begin{equation}\label{eq:kest}
P(\sigma^p>x)\sim c x^{-\alpha} ,\qquad \xto
\end{equation}
for some constant $c>0$. The index $\alpha$ is then obtained as the
unique positive solution to the equation $EA^\kappa=1$.
Relation \eqref{eq:kest} holds under
general conditions on $(A,B)$, which we do not give
here.
Regular variation of $(\sigma_t)$ is inherited by the solution to
\eqref{eq:5}.
\begin{lemma}\label{lem:kes}
Assume the conditions of Kesten \cite{kesten:1973} for the stochastic
recurrence equation \eqref{eq:5} and the moment conditions $EA^{\alpha+\vep}<\infty$ and
$EB^{\alpha+\vep}<\infty$ for some $\vep>0$. Then $(\sigma_t)$ is
regularly varying with index $\alpha p$ and for
${\boldsymbol\Sigma}_d=(\sigma_1,\ldots,\sigma_d)'$,
%
\begin{eqnarray}\label{eq:kestenlimit}
&&\frac{P(x^{-1}{\boldsymbol\Sigma}_d\in\cdot) }{P(\sigma>x)}\nonumber
\\[-8pt]
\\[-8pt]
&&\quad\stv
\alpha p\int_0^\infty t^{-\alpha p-1}
P\bigl(t(1,A_1^{1/p},\ldots,(A_{d-1}\cdots A_1)^{1/p})'\in
\cdot\bigr) \,\mathrm{d}t ,\qquad \xto.\nonumber
\end{eqnarray}
Moreover, if $E|Z|^{\alpha p+\delta}<\infty$ for some $\delta>0$, then
the stochastic volatility model $(X_t)$ is regularly varying with
index $\alpha p$, and the
limiting measure of
$\bfX_d=(X_1,\ldots,X_d)'$ is given by
%
\begin{eqnarray}\label{eq:kestenlimita}
&&\frac{P(x^{-1}\bfX_d\in\cdot) }{P(|X|>x)}\nonumber
\\[-8pt]
\\[-8pt]
&&\quad\stv
\frac{\alpha p}{ E|Z|^{\alpha p}}\int_0^\infty t^{-\alpha p-1}
P\bigl(t(Z_1,Z_2A_1^{1/p},\ldots,Z_d(A_{d-1}\cdots A_1)^{1/p})'\in
\cdot\bigr) \,\mathrm{d}t .\nonumber
\end{eqnarray}
If $Z$ is symmetric and $P(Z=0)=0$, then the limit
in \eqref{eq:kestenlimita} turns into
\[
\alpha p\int_0^\infty t^{-\alpha p-1}
P\bigl(t\bigl(\sign(Z_1),(Z_2/|Z_1|)A_1^{1/p},\ldots,(Z_d/|Z_1|)(A_{d-1}\cdots
A_1)^{1/p}\bigr)'\in
\cdot\bigr) \,\mathrm{d}t .
\]
\end{lemma}
\begin{pf} We take the approach in the proof of Corollary 2.7 in
\cite{basrak:davis:mikosch:2002}.
For every $t$, we have
%
\begin{equation}\label{eq:*}
\sigma_t^p = A_t\cdots A_1 \sigma_0^p+\sum_{i=1}^t A_t\cdots A_{i+1}
B_i ,
\end{equation}
and thus, applying the power operation component-wise,
\[
{\boldsymbol\Sigma}_d^p= \sigma_0^p (A_1,A_2A_1, \ldots,A_d\cdots
A_1)'+ \bfR_d ,
\]
where, by virtue of the moment conditions on $(A,B)$,
$E|\bfR_d|^{\alpha+\vep}<\infty$ for some $\vep>0$. By Kesten's
theorem (cf. \eqref{eq:kest}), as $\xto$,
\[
\frac{P(\sigma>xt)}{P(\sigma>x)}\to t^{-\alpha p}=\mu(t,\infty) ,\qquad
t>0 .
\]
Therefore, and in view of a version of the
multivariate Breiman result (see \cite{jessen:mikosch:2006}),
$P(x^{-1}{\boldsymbol\Sigma}_d\in\cdot)/P(\sigma>x)$ has
the same limit measure as
\begin{eqnarray*}
&&\frac{P(x^{-1}\sigma_0 (A_1,A_2A_1, \ldots,A_d\cdots A_1)'^{1/p}\in
\cdot)}{P(\sigma>x)}
\\
&&\quad\stv E\mu\{t>0\dvt t(A_1,A_2 A_1,\ldots,A_d\cdots
A_1)'^{1/p}\in\cdot\}
\\
&&\quad= \alpha p EA^\alpha\int_0^\infty t^{-\alpha p-1}
P\bigl(t(1,A_2,\ldots,A_d\cdots A_2)'^{1/p}\in\cdot\bigr) \,\mathrm{d}t
\\
&&\quad= \alpha p \int_0^\infty t^{-\alpha p-1}
P\bigl(t(1,A_1,\ldots,A_{d-1}\cdots A_1)'^{1/p}\in\cdot\bigr) \,\mathrm{d}t .
\end{eqnarray*}
Relation \eqref{eq:kestenlimita} follows by an application of the
multivariate Breiman result; compare Lemma~\ref{lem:1}.
\mbox{}\hfill\qed
\noqed\end{pf}
%
\subsection{Anticlustering condition}
\begin{lemma}
Assume that the conditions of Lemma~\ref{lem:kes} are satisfied,
ensuring that $(\sigma_t)$ is regularly varying with index $\alpha
p$.
Then the anticlustering condition \eqref{eq:ac} is satisfied for
$(\sigma_t)$ for a sequence $(r_n)$ satisfying
$r_n=\mathrm{O}(n^\gamma)$ for any small $\gamma>0$.
Moreover, if $E|Z|^{\alpha p+\delta}<\infty$ for some $\delta>0$,
then \eqref{eq:ac} also holds for $(X_t)$ with the same sequence $(r_n)$.
\end{lemma}
\begin{pf} Condition \eqref{eq:ac} for $(\sigma_t)$
follows from the proof of Theorem 2.10 in
\cite{basrak:davis:mikosch:2002}.
Indeed, \cite{basrak:davis:mikosch:2002} used
\eqref{eq:*} to show
that
%
\begin{equation}\label{eq:ll}
\lim_{m\to\infty}\limsup_{\nto} \sum_{m\le|t|\le r_n}
P(\sigma_t>a_n y\vert\sigma_0>a_n y)=0 ,\qquad y>0 .
\end{equation}
The corresponding result for $(X_t)$ follows along the lines of
the proof of \eqref{eq:ll}, exploiting \eqref{eq:*} and the
independence of
$(\sigma_t)$ and $(Z_t)$.
\end{pf}
%
\subsection{Main result for solution to stochastic recurrence equation}
We formulate an analog of Theorem~\ref{thm:mainsv}, summarizing the
foregoing results in the case of a solution to a stochastic recurrence equation.
\begin{theorem}\label{thm:sre}
Assume that the volatility sequence $(\sigma_t)$ is given via
the solution $(\sigma_t^p)$ of the
stochastic recurrence equation \eqref{eq:5} for some $p>0$. Assume
the following conditions:
\begin{itemize}
\item
$(\sigma_t^p)$ satisfies Kesten's \cite{kesten:1973} conditions.
\item
$(\sigma_t)$ is strongly mixing with
geometric rate.
\end{itemize}
Then
\begin{enumerate}[(3)]
\item[(1)]
$(\sigma_t)$ is
regularly varying with index $\alpha p$ and limiting measures given in
\eqref{eq:kestenlimit}.\eject
\item[(2)]
Condition ${\mathcal A}(a_n)$ is satisfied for any
$(r_n)$ satisfying $r_n=\mathrm{o}(n)$ and $r_n\ge c\log n$ for some $c>0$.
\item[(3)]
The anticlustering
condition \eqref{eq:ac} holds for $(\sigma_t)$ with a sequence $(r_n)$
satisfying $r_n=\mathrm{O}(n^\gamma)$ for any small $\gamma>0$.
\end{enumerate}
Moreover, if
$E|Z|^{\alpha p+\delta}<\infty$ for some $\delta>0$, then the
following hold:
\begin{enumerate}[(6)]
\item[(4)]
$(X_t)$ is
regularly varying with index $\alpha p$ and limiting measures given in
\eqref{eq:kestenlimita}.
\item[(5)]
$(X_t)$ is strongly mixing with
geometric rate, and condition ${\mathcal A}(a_n)$ is satisfied for any
$(r_n)$ satisfying $r_n=\mathrm{o}(n)$ and $r_n\ge c\log n$ for some $c>0$.
\item[(6)]
The anticlustering
condition \eqref{eq:ac} holds for $(X_t)$
and sequences $(r_n)$
satisfying $r_n=\mathrm{O}(n^\gamma)$ for any small $\gamma>0$.
\end{enumerate}
In particular, Theorem~\ref{thm:dh} is applicable to the sequences
$(\sigma_t)$ and $(X_t)$.
\end{theorem}

Now we can again use the machinery
of
\cite
{davis:hsing:1995,davis:mikosch:1998,davis:mikosch:2001,basrak:davis:mikosch:2002}
to derive various limit results for functionals of the sequence $(X_t)$.
We only derive the extremal index of $(X_t)$ in a special
situation, to show the
crucial difference between the exponential $\operatorname{AR}(1)$ process considered
in Section~\ref{subsec:ar} and the present situation.
%

\begin{figure}[b]

\includegraphics{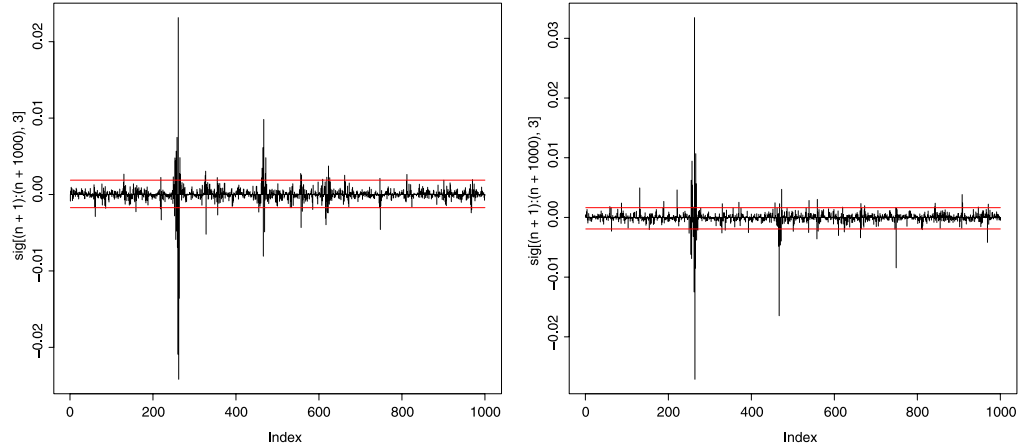}

  \caption{\emph{Left:} 1000 realizations of a \garch\ process with
parameters
$\alpha_0=10^{-7}$, $\alpha_1=0.1$, $\beta_1=0.89$ and i.i.d. standard
normal noise.
\emph{Right:} Realizations of a stochastic volatility model, where $(\sigma_t)$ is taken from the \garch\ process in the
left graph and $Z$ is
standard normal.
In both graphs, $(\sigma_t)$ and $(X_t)$ are regularly varying with index 4 causing extremal clustering in both sequences.
The parallel lines indicate the 0.99 and 0.01 quantiles of the
distribution of $X$.}\label{fig:2}
\end{figure}

\begin{example}
Recall the definition of a \garch\ process from
Example~\ref{exam:garch}.
We assume that
$(\sigma_t^2)$ is the squared volatility process of a \garch\
process that is regularly varying with index $\alpha>0$ and
$E|Z|^{2\alpha+\delta}<\infty$ for some $\delta>0$. Such a \garch\
process and the corresponding
stochastic volatility model are shown in Figure~\ref{fig:2}.
Assume that both $(\sigma_t)$ and the corresponding stochastic
volatility model
$(X_t)$ satisfy the conditions of Theorem~\ref{thm:dh}; sufficient
conditions are given in Theorem~\ref{thm:sre}. It is well known
(e.g., \cite{mikosch:starica:2000}, Theorem 4.1) that
the extremal indices of $(\sigma_t)$ and $(\sigma_t^2)$ coincide and are
given by
\[
\theta_\sigma= \alpha\int_1^{\infty} P \Biggl(\sup_{t\ge1} \prod
_{j=1}^t A_j\le y^{-1} \Biggr) y^{-1-\alpha} \,\mathrm{d}y ,
\]
where $A_j=\alpha_1\eta_j^2+\beta_1$, $j\ge1$.
For the extremal index $\theta_{|X|}$ of the sequences $(|X_t|)$ and
$(X_t^2)$, we use the expression in \cite{davis:hsing:1995} given by
%
\begin{equation}\label{eq:index}
\theta_{|X|}=\lim_{m\to\infty}
\frac{
(|\theta_0^{(m)}|^{\alpha}-\max_{j=1,\ldots,m}
|\theta_j^{(m)}|^{\alpha} )_+}{E|\theta_0^{(m)}|^{\alpha}} ,
\end{equation}
where $\bfTh^{(m)}=(\theta_j^{(m)})_{|j|\le m}$ is a vector with
values in the
unit sphere $\bbs^{2m}$ of $\bbr^{2m+1}$, which has the spectral distribution
 of the random vector $\wh\bfX^{(m)}=(X_t^2)_{|t|\le m}$, that is,
\[
\frac{P(|\wh\bfX^{(m)}|>x,\wh\bfX^{(m)}/|\wh\bfX^{(m)}|\in\cdot)}
{P(|\wh\bfX^{(m)}|>x)}\stw P\bigl(\bfTh^{(m)}\in\cdot\bigr) , \qquad\xto.
\]
For any Borel set $S\subset\bbs^{2m}$ that is a continuity set with
respect to $P(\bfTh^{(m)}\in\cdot)$, we conclude from \eqref{eq:kestenlimita}
with $\bfR^{(m)}=(Z_1^2,Z_2^2A_1,\ldots,Z_{2m+1} ^2A_{2m}\cdots A_1)'$
\begin{eqnarray*}
\frac{
P(|\wh\bfX^{(m)}|>x,\wh\bfX^{(m)}/|\wh\bfX^{(m)}|\in S )}
{P(|\wh\bfX^{(m)}|>x)}
&\to&
\frac{\alpha\int_0^\infty t^{-\alpha-1}
P(t |\bfR^{(m)}| I_{\{
\bfR^{(m)}/|\bfR^{(m)}|\in S\}}>1) \,\mathrm{d}t}{E|\bfR^{(m)}|^{\alpha}}
\\
&=& \frac{E|\bfR^{(m)}|^{\alpha} I_{\{
\bfR^{(m)}/|\bfR^{(m)}|\in S\}}
}{E|\bfR^{(m)}|^{\alpha}}=P\bigl(\bfTh^{(m)}\in S\bigr) .
\end{eqnarray*}
The latter relation, \eqref{eq:index}, and the fact that $EA^\alpha=1$
yield
\[
\theta_{|X|}=\lim_{m\to\infty}\frac{E (
|Z_1|^{2\alpha}-\max_{j=2,\ldots,m} (Z_j^2 \prod_{i=2}^j
A_i)^\alpha)_+}{E|Z|^{2\alpha}} .
\]
A comparison with Theorem~4.1 in \cite{mikosch:starica:2000} shows
that a similar expression can be derived for the extremal index
$\theta_{|X|}$ of the \garch\ process; the
$Z$'s must be replaced by the corresponding $\eta$'s.
(For details on the foregoing calculations, see
\cite{mikosch:starica:2000}.)
A direct comparison of the magnitude of the extremal indices of a
\garch\ process and the corresponding stochastic volatility model
seems difficult.
\end{example}
%
\section{Moving average processes}\label{subsec:3}
In this section we assume that the
volatility process $(\sigma_t)$ is given in the form $\sigma
_t^p=|Y_t|$ for some
$p>0$ and $Y_t=\sum_{j=0}^q \psi_j \eta_{t-j}$, $t\in\bbz$, for
some $q\ge
1$ and an i.i.d. sequence $(\eta_t)$ such that $\eta$ is regularly varying in the sense of
\eqref{eq:rv} with tail balance coefficients $\wt p,\wt q\ge0$, $\wt
p+\wt q=1$ and index $\alpha>0$.
Because $\bfY_d=(Y_1,\ldots,Y_d)'$ has representation as a
linear transformation of a finite
vector of the $Z$'s, an application of the continuous mapping theorem
implies that the vector $\bfY_d$ is regularly varying with index
$\alpha$. Writing $\psi_d=0$ for $d\notin\{0,\ldots,q\}$, we
conclude from \cite{davis:resnick:1985}, Theorem 2.4, that
%
\begin{eqnarray}\label{eq:dr}
\frac{P(x^{-1}\bfY_d \in\cdot)}{P(|\eta|>x)}&\stv&
\alpha\sum_{j=0}^{q+d-1}\int_{\ov\bbr_0}
|x|^{-\alpha-1} \bigl[\wt p I_{(0,\infty)}(x)
+\wt q I_{(-\infty,0)}(x)\bigr]\nonumber
\\[-8pt]
\\[-8pt]
&&\hphantom{\alpha\sum_{j=0}^{q+d-1}\int_{\ov\bbr_0}}{}\times I_{\{x
(\psi_{j-d+1},\ldots,\psi_{j})\in\cdot\}} \,\mathrm{d}x .\nonumber
\end{eqnarray}
The mixing condition ${\mathcal A}(a_n)$ and the anticlustering
condition \eqref{eq:ac} are automatically satisfied for $(\sigma_t)$
and $(X_t)$. Thus
Theorem~\ref{thm:main} holds. We conclude from \eqref{eq:dr} and
Breiman's result that
\begin{eqnarray*}
&&\frac{P(x^{-1}(Z_1\sigma_1,\ldots,Z_d\sigma_d)\in\cdot)}{P(|\eta
|>x)}
\\
&&\quad\stv
\alpha p \sum_{j=0}^{q+d-1}\int_0^\infty
|x|^{-\alpha p-1} P\bigl(x
(Z_1|\psi_{j-d+1}|^{1/p},\ldots,Z_d|\psi_{j}|^{1/p})\in\cdot\bigr)\,\mathrm{d}x .
\end{eqnarray*}
An application of \eqref{eq:index} yields
%
\begin{equation}\label{eq:li}
\theta_{|X|}=\frac{E\max_{j=0,\ldots,q} |Z_j|^{\alpha p}
|\psi_j|^\alpha}{E|Z|^{\alpha p}\sum_{j=0}^q|\psi_j|^\alpha} .
\end{equation}
In the degenerate case when $Z=1$, we get the well-known form of
the extremal index of the absolute values of a moving average process
(see \cite{davis:resnick:1985};
cf. \cite{embrechts:kluppelberg:mikosch:1997},
page 415). Again, a direct comparison of the value \eqref{eq:li} with the
corresponding one for $Z=1$ seems difficult.

\begin{remark}
The foregoing techniques can be
applied in the case where $(Y_t)$ constitutes an infinite moving
average process as well. However, in this case mixing conditions are
generally difficult to check; instead, \cite{davis:resnick:1985}
used approximations of an infinite moving average by finite moving
averages. This technique does not completely fit into the framework of
\cite{davis:hsing:1995}; see Theorem~\ref{thm:dh} above. However, if $(Y_t)$
is an ARMA process with i.i.d. noise
$(\eta_t)$ that is regularly varying with index $\alpha>1$ and has
a positive density in some neighborhood of $E\eta$, then $(Y_t)$
is strongly mixing with geometric rate. Then ${\mathcal A}(a_n)$
holds for every sequence $(r_n)$ with $r\ge c\log n$ for some $c>0$, and
${\mathcal A}(a_n)$ also holds for $(X_t)$ and the same sequence $(r_n)$. The anticlustering condition for $(X_t)$ can be checked in
this case as well, but the calculations are lengthy. We omit further details.
\end{remark}
%
\section{Concluding remarks}
The aim of this paper was to show that the stochastic volatility model
$(X_t)$ given by
\eqref{eq:sv} may exhibit extremal clustering provided that $(\sigma
_t)$ is
a regularly varying sequence with index $\alpha>0$ and the i.i.d.
noise sequence $(Z_t)$ has $(\alpha+\vep)$th moment for some $\vep>0$. Extremal
clustering is inherited from the volatility sequence $(\sigma_t)$. If
$(\sigma_t)$ does not have extremal clusters, then neither does the
sequence $(X_t)$. An example of this lack of clustering is given by an
exponential $\operatorname{AR}(1)$ process $\sigma_t=\mathrm{e}^{Y_t}$, $Y_t=\varphi
Y_{t-1}+\eta_t$ for $\varphi\in(-1,1)$ and an i.i.d. regularly
varying sequence
 $(\mathrm{e}^{\eta_t})$. The results of Section~\ref{subsec:ar} show that
the sequence $(X_t)$ above high levels essentially behaves like the
i.i.d. sequence $(\mathrm{e}^{\eta_t})$, resulting in an extremal index
$\theta_{|X|}=1$.
This is surprising, given that the autocorrelation
function of $(|X_t|)$ is not negligible. This example includes
$(\sigma_t)$ given by the dynamics of an EGARCH process. The EGARCH
process itself does then not exhibit extremal clustering either.

In contrast to an exponential $\operatorname{AR}(1)$, the stochastic volatility model
\eqref{eq:sv}
exhibits extremal clustering if the dynamics of
$(\sigma_t)$ or some positive power of
it are given by a moving average or the solution to a stochastic
recurrence equation. The latter
case captures
the example of the volatility sequence of a \garch\ process.

We have chosen to describe extremal clustering in terms of the
extremal index of the sequence $(X_t)$. If $\theta_{|X|}<1$, then
evaluating this quantity is difficult in the examples considered. We
would depend on
numerical or Monte Carlo methods if we were interested in numerical
values of $\theta_{|X|}$. These methods also would depend on the model.

The literature on the extremes of the
stochastic volatility model focuses on the case where $(\sigma_t)$ is
lognormal and
$(Z_t)$ is i.i.d. normal or regularly varying (cf. \cite{davis:mikosch:2009}).
In these cases, $(X_t)$ does not have extremal clusters. The latter
property can be considered a disadvantage for modeling return
series that are believed to have the clustering property, often
referred to as volatility clusters. From a modeling standpoint,
neither the stochastic volatility model with or without extremal
clusters nor any
standard model such as GARCH or EGARCH can be discarded as long as
no efficient methods for distinguishing between these models exist.
For example, the volatility dynamics of an EGARCH model and a
stochastic volatility model
with exponential $\operatorname{AR}(1)$ volatility
are rather similar and so are the volatility dynamics of
a \garch\ and a stochastic volatility model with \garch\ volatility.

\section*{Acknowledgements}
 This paper was written
when Mohsen Rezapour visited the
Department of Mathematics at the University of Copenhagen
November 2010--April 2011.
He would like to thank the
Department of Mathematics for hospitality.
He also thanks the Office of Graduate Studies at the University of
Isfahan for its support. Thomas Mikosch's research
is partly supported by the Danish
Natural Science Research Council (FNU) Grants
09-072331 ``Point process modelling and statistical inference''
and 10-084172 ``Heavy tail phenomena: Modeling and estimation.''
The constructive remarks of the referees and Associate Editor led to
an improved presentation of the paper. We would like to thank
them.\looseness=1
%

\printhistory


\begin{thebibliography}{32}

\bibitem{andersen:davis:kreiss:mikosch:2009}
\begin{bbook}[auto:STB|2012/03/21|07:41:58]
\bauthor{\bsnm{Andersen},~\bfnm{T.~G.}\binits{T.G.}},
  \bauthor{\bsnm{Davis},~\bfnm{R.~A.}\binits{R.A.}},
  \bauthor{\bsnm{Kreiss},~\bfnm{J.~P.}\binits{J.P.}} \AND
  \bauthor{\bsnm{Mikosch},~\bfnm{T.}\binits{T.}}
(\byear{2009}).
\btitle{Handbook of Financial Time Series}.
\baddress{Berlin}: \bpublisher{Springer}.
\bptok{imsref}%
\end{bbook}
\endbibitem

\bibitem{basrak:davis:mikosch:2002a}
\begin{barticle}[mr]
\bauthor{\bsnm{Basrak},~\bfnm{Bojan}\binits{B.}},
  \bauthor{\bsnm{Davis},~\bfnm{Richard~A.}\binits{R.A.}} \AND
  \bauthor{\bsnm{Mikosch},~\bfnm{Thomas}\binits{T.}}
(\byear{2002}).
\btitle{A characterization of multivariate regular variation}.
\bjournal{Ann. Appl. Probab.}
\bvolume{12}
\bpages{908--920}.
\bid{doi={10.1214/aoap/1031863174}, issn={1050-5164}, mr={1925445}}
\bptok{imsref}%
\end{barticle}
\endbibitem

\bibitem{basrak:davis:mikosch:2002}
\begin{barticle}[mr]
\bauthor{\bsnm{Basrak},~\bfnm{Bojan}\binits{B.}},
  \bauthor{\bsnm{Davis},~\bfnm{Richard~A.}\binits{R.A.}} \AND
  \bauthor{\bsnm{Mikosch},~\bfnm{Thomas}\binits{T.}}
(\byear{2002}).
\btitle{Regular variation of {GARCH} processes}.
\bjournal{Stochastic Process. Appl.}
\bvolume{99}
\bpages{95--115}.
\bid{doi={10.1016/S0304-4149(01)00156-9}, issn={0304-4149}, mr={1894253}}
\bptok{imsref}%
\end{barticle}
\endbibitem

\bibitem{basrak:segers:2009}
\begin{barticle}[mr]
\bauthor{\bsnm{Basrak},~\bfnm{Bojan}\binits{B.}} \AND
  \bauthor{\bsnm{Segers},~\bfnm{Johan}\binits{J.}}
(\byear{2009}).
\btitle{Regularly varying multivariate time series}.
\bjournal{Stochastic Process. Appl.}
\bvolume{119}
\bpages{1055--1080}.
\bid{doi={10.1016/j.spa.2008.05.004}, issn={0304-4149}, mr={2508565}}
\bptnote{check year}
\bptok{imsref}%
\end{barticle}
\endbibitem

\bibitem{bingham:goldie:teugels:1987}
\begin{bbook}[mr]
\bauthor{\bsnm{Bingham},~\bfnm{N.~H.}\binits{N.H.}},
  \bauthor{\bsnm{Goldie},~\bfnm{C.~M.}\binits{C.M.}} \AND
  \bauthor{\bsnm{Teugels},~\bfnm{J.~L.}\binits{J.L.}}
(\byear{1987}).
\btitle{Regular Variation}.
\bseries{Encyclopedia of Mathematics and Its Applications}
\bvolume{27}.
\baddress{Cambridge}: \bpublisher{Cambridge Univ. Press}.
\bid{mr={0898871}}
\bptok{imsref}%
\end{bbook}
\endbibitem

\bibitem{boman:lindskog:2007}
\begin{barticle}[mr]
\bauthor{\bsnm{Boman},~\bfnm{Jan}\binits{J.}} \AND
  \bauthor{\bsnm{Lindskog},~\bfnm{Filip}\binits{F.}}
(\byear{2009}).
\btitle{Support theorems for the {R}adon transform and {C}ram\'er-{W}old
  theorems}.
\bjournal{J. Theoret. Probab.}
\bvolume{22}
\bpages{683--710}.
\bid{doi={10.1007/s10959-008-0151-0}, issn={0894-9840}, mr={2530109}}
\bptnote{check year}
\bptok{imsref}%
\end{barticle}
\endbibitem

\bibitem{breidt:davis:1998}
\begin{barticle}[mr]
\bauthor{\bsnm{Breidt},~\bfnm{F.~Jay}\binits{F.J.}} \AND
  \bauthor{\bsnm{Davis},~\bfnm{Richard~A.}\binits{R.A.}}
(\byear{1998}).
\btitle{Extremes of stochastic volatility models}.
\bjournal{Ann. Appl. Probab.}
\bvolume{8}
\bpages{664--675}.
\bid{doi={10.1214/aoap/1028903446}, issn={1050-5164}, mr={1627756}}
\bptok{imsref}%
\end{barticle}
\endbibitem

\bibitem{breiman:1965}
\begin{barticle}[auto:STB|2012/03/21|07:41:58]
\bauthor{\bsnm{Breiman},~\bfnm{L.}\binits{L.}}
(\byear{1965}).
\btitle{On some limit theorems similar to the arc-sin law}.
\bjournal{Theory Probab. Appl.}
\bvolume{10}
\bpages{323--331}.
\bptok{imsref}%
\end{barticle}
\endbibitem

\bibitem{davis:resnick:1985}
\begin{barticle}[mr]
\bauthor{\bsnm{Davis},~\bfnm{Richard~A.}\binits{R.A.}} \AND
  \bauthor{\bsnm{Resnick},~\bfnm{Sidney}\binits{S.}}
(\byear{1985}).
\btitle{Limit theory for moving averages of random variables with regularly
  varying tail probabilities}.
\bjournal{Ann. Probab.}
\bvolume{13}
\bpages{179--195}.
\bid{issn={0091-1798}, mr={0770636}}
\bptok{imsref}%
\end{barticle}
\endbibitem

\bibitem{davis:hsing:1995}
\begin{barticle}[mr]
\bauthor{\bsnm{Davis},~\bfnm{Richard~A.}\binits{R.A.}} \AND
  \bauthor{\bsnm{Hsing},~\bfnm{Tailen}\binits{T.}}
(\byear{1995}).
\btitle{Point process and partial sum convergence for weakly dependent random
  variables with infinite variance}.
\bjournal{Ann. Probab.}
\bvolume{23}
\bpages{879--917}.
\bid{issn={0091-1798}, mr={1334176}}
\bptok{imsref}%
\end{barticle}
\endbibitem

\bibitem{davis:mikosch:1998}
\begin{barticle}[mr]
\bauthor{\bsnm{Davis},~\bfnm{Richard~A.}\binits{R.A.}} \AND
  \bauthor{\bsnm{Mikosch},~\bfnm{Thomas}\binits{T.}}
(\byear{1998}).
\btitle{The sample autocorrelations of heavy-tailed processes with applications
  to {ARCH}}.
\bjournal{Ann. Statist.}
\bvolume{26}
\bpages{2049--2080}.
\bid{doi={10.1214/aos/1024691368}, issn={0090-5364}, mr={1673289}}
\bptok{imsref}%
\end{barticle}
\endbibitem

\bibitem{davis:mikosch:2001}
\begin{barticle}[mr]
\bauthor{\bsnm{Davis},~\bfnm{Richard~A.}\binits{R.A.}} \AND
  \bauthor{\bsnm{Mikosch},~\bfnm{Thomas}\binits{T.}}
(\byear{2001}).
\btitle{Point process convergence of stochastic volatility processes with
  application to sample autocorrelation}.
\bjournal{J. Appl. Probab.}
\bvolume{38A}
\bpages{93--104}.
\bid{doi={10.1239/jap/1085496594}, issn={0021-9002}, mr={1915537}}
\bptok{imsref}%
\end{barticle}
\endbibitem

\bibitem{davis:mikosch:2009}
\begin{bincollection}[auto:STB|2012/03/21|07:41:58]
\bauthor{\bsnm{Davis},~\bfnm{R.~A.}\binits{R.A.}} \AND
  \bauthor{\bsnm{Mikosch},~\bfnm{T.}\binits{T.}}
(\byear{2009}).
\btitle{Extremes of stochastic volatility models}.
In \bbooktitle{Handbook of Financial Time Series}
(\beditor{\bfnm{T.~G.}\binits{T.G.}~\bsnm{Andersen}},
\beditor{\bfnm{R.~A.}\binits{R.A.}~\bsnm{Davis}},
\beditor{\bfnm{J.~P.}\binits{J.P.}~\bsnm{Kreiss}} \AND
   \beditor{\bfnm{T.}\binits{T.}~\bsnm{Mikosch}}, eds.)
\bpages{355--364}.
\baddress{Berlin}: \bpublisher{Springer}.
\bptok{imsref}%
\end{bincollection}
\endbibitem

\bibitem{davis:mikosch:2009a}
\begin{bincollection}[auto:STB|2012/03/21|07:41:58]
\bauthor{\bsnm{Davis},~\bfnm{R.~A.}\binits{R.A.}} \AND
  \bauthor{\bsnm{Mikosch},~\bfnm{T.}\binits{T.}}
(\byear{2009}).
\btitle{Fundamental properties of stochastic volatility models}.
In \bbooktitle{Handbook of Financial Time Series}
(\beditor{\bfnm{T.~G.}\binits{T.G.}~\bsnm{Andersen}},
  \beditor{\bfnm{R.~A.}\binits{R.A.}~\bsnm{Davis}},
  \beditor{\bfnm{J.~P.}\binits{J.P.}~\bsnm{Kreiss}} \AND
  \beditor{\bfnm{T.}\binits{T.}~\bsnm{Mikosch}}, eds.)
\bpages{255--267}.
\baddress{Berlin}: \bpublisher{Springer}.
\bptok{imsref}%
\end{bincollection}
\endbibitem

\bibitem{doukhan:1994}
\begin{bbook}[mr]
\bauthor{\bsnm{Doukhan},~\bfnm{Paul}\binits{P.}}
(\byear{1994}).
\btitle{Mixing: Properties and Examples}.
\bseries{Lecture Notes in Statistics}
\bvolume{85}.
\baddress{New York}: \bpublisher{Springer}.
\bid{mr={1312160}}
\bptok{imsref}%
\end{bbook}
\endbibitem

\bibitem{embrechts:kluppelberg:mikosch:1997}
\begin{bbook}[mr]
\bauthor{\bsnm{Embrechts},~\bfnm{Paul}\binits{P.}},
  \bauthor{\bsnm{Kl{\"u}ppelberg},~\bfnm{Claudia}\binits{C.}} \AND
  \bauthor{\bsnm{Mikosch},~\bfnm{Thomas}\binits{T.}}
(\byear{1997}).
\btitle{Modelling Extremal Events: For Insurance and Finance}.
\bseries{Applications of Mathematics (New York)}
\bvolume{33}.
\baddress{Berlin}: \bpublisher{Springer}.
\bid{mr={1458613}}
\bptok{imsref}%
\end{bbook}
\endbibitem

\bibitem{embrechts:veraverbeke:1982}
\begin{barticle}[mr]
\bauthor{\bsnm{Embrechts},~\bfnm{P.}\binits{P.}} \AND
  \bauthor{\bsnm{Veraverbeke},~\bfnm{N.}\binits{N.}}
(\byear{1982}).
\btitle{Estimates for the probability of ruin with special emphasis on the
  possibility of large claims}.
\bjournal{Insurance Math. Econom.}
\bvolume{1}
\bpages{55--72}.
\bid{doi={10.1016/0167-6687(82)90021-X}, issn={0167-6687}, mr={0652832}}
\bptok{imsref}%
\end{barticle}
\endbibitem

\bibitem{goldie:1991}
\begin{barticle}[mr]
\bauthor{\bsnm{Goldie},~\bfnm{Charles~M.}\binits{C.M.}}
(\byear{1991}).
\btitle{Implicit renewal theory and tails of solutions of random equations}.
\bjournal{Ann. Appl. Probab.}
\bvolume{1}
\bpages{126--166}.
\bid{issn={1050-5164}, mr={1097468}}
\bptok{imsref}%
\end{barticle}
\endbibitem

\bibitem{hult:lindskog:2005}
\begin{barticle}[mr]
\bauthor{\bsnm{Hult},~\bfnm{Henrik}\binits{H.}} \AND
  \bauthor{\bsnm{Lindskog},~\bfnm{Filip}\binits{F.}}
(\byear{2005}).
\btitle{Extremal behavior of regularly varying stochastic processes}.
\bjournal{Stochastic Process. Appl.}
\bvolume{115}
\bpages{249--274}.
\bid{doi={10.1016/j.spa.2004.09.003}, issn={0304-4149}, mr={2111194}}
\bptok{imsref}%
\end{barticle}
\endbibitem

\bibitem{hult:lindskog:2006}
\begin{barticle}[mr]
\bauthor{\bsnm{Hult},~\bfnm{Henrik}\binits{H.}} \AND
  \bauthor{\bsnm{Lindskog},~\bfnm{Filip}\binits{F.}}
(\byear{2006}).
\btitle{On {K}esten's counterexample to the {C}ram\'er--{W}old device for
  regular variation}.
\bjournal{Bernoulli}
\bvolume{12}
\bpages{133--142}.
\bid{issn={1350-7265}, mr={2202325}}
\bptok{imsref}%
\end{barticle}
\endbibitem

\bibitem{hult:lindskog:2006a}
\begin{barticle}[mr]
\bauthor{\bsnm{Hult},~\bfnm{Henrik}\binits{H.}} \AND
  \bauthor{\bsnm{Lindskog},~\bfnm{Filip}\binits{F.}}
(\byear{2006}).
\btitle{Regular variation for measures on metric spaces}.
\bjournal{Publ. Inst. Math. (Beograd) (N.S.)}
\bvolume{80(94)}
\bpages{121--140}.
\bid{doi={10.2298/PIM0694121H}, issn={0350-1302}, mr={2281910}}
\bptok{imsref}%
\end{barticle}
\endbibitem

\bibitem{jessen:mikosch:2006}
\begin{barticle}[mr]
\bauthor{\bsnm{Jessen},~\bfnm{Anders~Hedegaard}\binits{A.H.}} \AND
  \bauthor{\bsnm{Mikosch},~\bfnm{Thomas}\binits{T.}}
(\byear{2006}).
\btitle{Regularly varying functions}.
\bjournal{Publ. Inst. Math. (Beograd) (N.S.)}
\bvolume{80(94)}
\bpages{171--192}.
\bid{doi={10.2298/PIM0694171J}, issn={0350-1302}, mr={2281913}}
\bptok{imsref}%
\end{barticle}
\endbibitem

\bibitem{kesten:1973}
\begin{barticle}[mr]
\bauthor{\bsnm{Kesten},~\bfnm{Harry}\binits{H.}}
(\byear{1973}).
\btitle{Random difference equations and renewal theory for products of random
  matrices}.
\bjournal{Acta Math.}
\bvolume{131}
\bpages{207--248}.
\bid{issn={0001-5962}, mr={0440724}}
\bptok{imsref}%
\end{barticle}
\endbibitem

\bibitem{kluppelberg:pergam:2007}
\begin{barticle}[mr]
\bauthor{\bsnm{Kl{\"u}ppelberg},~\bfnm{Claudia}\binits{C.}} \AND
  \bauthor{\bsnm{Pergamenchtchikov},~\bfnm{Serguei}\binits{S.}}
(\byear{2007}).
\btitle{Extremal behaviour of models with multivariate random recurrence
  representation}.
\bjournal{Stochastic Process. Appl.}
\bvolume{117}
\bpages{432--456}.
\bid{doi={10.1016/j.spa.2006.09.001}, issn={0304-4149}, mr={2305380}}
\bptok{imsref}%
\end{barticle}
\endbibitem

\bibitem{kulik:soulier:2011}
\begin{barticle}[mr]
\bauthor{\bsnm{Kulik},~\bfnm{Rafa{\l}}\binits{R.}} \AND
  \bauthor{\bsnm{Soulier},~\bfnm{Philippe}\binits{P.}}
(\byear{2011}).
\btitle{The tail empirical process for long memory stochastic volatility
  sequences}.
\bjournal{Stochastic Process. Appl.}
\bvolume{121}
\bpages{109--134}.
\bid{doi={10.1016/j.spa.2010.09.001}, issn={0304-4149}, mr={2739008}}
\bptok{imsref}%
\end{barticle}
\endbibitem

\bibitem{leadbetter:lindgren:rootzen:1983}
\begin{bbook}[mr]
\bauthor{\bsnm{Leadbetter},~\bfnm{M.~R.}\binits{M.R.}},
  \bauthor{\bsnm{Lindgren},~\bfnm{Georg}\binits{G.}} \AND
  \bauthor{\bsnm{Rootz{\'e}n},~\bfnm{Holger}\binits{H.}}
(\byear{1983}).
\btitle{Extremes and Related Properties of Random Sequences and Processes}.
\bseries{Springer Series in Statistics}.
\baddress{New York}: \bpublisher{Springer}.
\bid{mr={0691492}}
\bptok{imsref}%
\end{bbook}
\endbibitem

\bibitem{mikosch:starica:2000}
\begin{barticle}[mr]
\bauthor{\bsnm{Mikosch},~\bfnm{Thomas}\binits{T.}} \AND
  \bauthor{\bsnm{St{\u{a}}ric{\u{a}}},~\bfnm{C{\u{a}}t{\u{a}}lin}\binits{C.}}
(\byear{2000}).
\btitle{Limit theory for the sample autocorrelations and extremes of a
  {$\operatorname{GARCH}(1,1)$} process}.
\bjournal{Ann. Statist.}
\bvolume{28}
\bpages{1427--1451}.
\bid{doi={10.1214/aos/1015957401}, issn={0090-5364}, mr={1805791}}
\bptok{imsref}%
\end{barticle}
\endbibitem

\bibitem{mokkadem:1990}
\begin{barticle}[mr]
\bauthor{\bsnm{Mokkadem},~\bfnm{Abdelkader}\binits{A.}}
(\byear{1990}).
\btitle{Propri\'et\'es de m\'elange des processus autor\'egressifs
  polynomiaux}.
\bjournal{Ann. Inst. Henri Poincar\'e Probab. Stat.}
\bvolume{26}
\bpages{219--260}.
\bid{issn={0246-0203}, mr={1063750}}
\bptok{imsref}%
\end{barticle}
\endbibitem

\bibitem{nelson:1991}
\begin{barticle}[mr]
\bauthor{\bsnm{Nelson},~\bfnm{Daniel~B.}\binits{D.B.}}
(\byear{1991}).
\btitle{Conditional heteroskedasticity in asset returns: A new approach}.
\bjournal{Econometrica}
\bvolume{59}
\bpages{347--370}.
\bid{doi={10.2307/2938260}, issn={0012-9682}, mr={1097532}}
\bptok{imsref}%
\end{barticle}
\endbibitem

\bibitem{resnick:1987}
\begin{bbook}[mr]
\bauthor{\bsnm{Resnick},~\bfnm{Sidney~I.}\binits{S.I.}}
(\byear{1987}).
\btitle{Extreme Values, Regular Variation, and Point Processes}.
\bseries{Applied Probability. A~Series of the Applied Probability Trust}
\bvolume{4}.
\baddress{New York}: \bpublisher{Springer}.
\bid{mr={0900810}}
\bptok{imsref}%
\end{bbook}
\endbibitem

\bibitem{resnick:2007}
\begin{bbook}[mr]
\bauthor{\bsnm{Resnick},~\bfnm{Sidney~I.}\binits{S.I.}}
(\byear{2007}).
\btitle{Heavy-Tail Phenomena: Probabilistic and Statistical Modeling}.
\bseries{Springer Series in Operations Research and Financial Engineering}.
\baddress{New York}: \bpublisher{Springer}.
\bid{mr={2271424}}
\bptok{imsref}%
\end{bbook}
\endbibitem

\end{thebibliography}
\end{document}